%% file: CPRS-2026-Jul.tex
\documentclass[reqno]{amsart} % reqno coloca as equações na direita
\usepackage{comment}
\usepackage{amsfonts,color}
\usepackage{latexsym,amssymb}
\usepackage{amsmath}
\usepackage[utf8]{inputenc}
\usepackage{enumitem}
\usepackage{xcolor}
\usepackage{fix-cm}

\makeatletter
\renewcommand{\tagform@}[1]{\maketag@@@{\color{blue}(#1)}}
\makeatother

\newcommand{\Oe}{\Omega_{\ell}}
\newcommand{\Onl}{\Omega_{n\ell}}
\newcommand{\Lop}{\mathcal{L}}

\newcommand{\normaind}[2]{\left\lVert #1 \right\rVert_{#2}}
\newtheorem{theorem}{Theorem}[section]
\newtheorem{lemma}[theorem]{Lemma}
\newtheorem{proposition}[theorem]{Proposition}
\newtheorem{corollary}[theorem]{Corollary}
\newtheorem{remark}[theorem]{Remark}
\newtheorem{definition}[theorem]{Definition}

\usepackage[bookmarksnumbered,colorlinks, linkcolor=blue, citecolor=red, pagebackref, bookmarks, breaklinks]{hyperref}
\usepackage{tikz-cd} %diagramas
\usepackage{xparse,tikz}
\usetikzlibrary{positioning,calc}
\usetikzlibrary {shapes.geometric}
\usetikzlibrary {patterns}
\usetikzlibrary {patterns,patterns.meta}
\usepackage[hyperpageref]{backref}

\catcode`@=11 \@addtoreset{equation}{section} \catcode`@=12
	
\begin{document}

\title[Maximum Principles and Principal Eigenvalue]{Maximum Principles and Principal Eigenvalue for Coupled Local/Nonlocal Models}

\author[W. Cintra]{Willian Cintra$^{1}$}
\address{$^{1}$Universidade de Brasília, Departamento de Matem\'atica, BR-70910900 Brasília-DF, Brazil.  E-mails: willian@unb.br, jonatasperalta58@gmail.com}

	\author[J.  Peralta]{Jônatas Peralta$^1$}
	
    \author[J. Rossi]{Julio D. Rossi$^2$}
    \address{$^2$Departamento de Matemáticas y Estadística,  Universidad  Torcuato Di Tella, Buenos Aires, Argentina.
    E-mail: julio.rossi@utdt.edu}
    
    \author[A. Suárez]{Antonio Suárez$^{3}$}
\address{$^{3}$Departamento de Ecuaciones Diferenciales y Análisis Numérico and IMUS, Univ. de Sevilla, Avd. Reina Mercedes, s/n, 41013, Sevilla, Spain. E-mail: suarez@us.es}
		
\begin{abstract} 
In this paper, we analyze the validity of the maximum principle and the existence of a principal eigenvalue for a problem involving two regions, in one of them the main operator is the classical Laplacian (a local operator) while in the other region the operator is a nonlocal one (these regions are coupled  via a transmission term that is also nonlocal). The main difficulty here comes from the lack  of compactness of the nonlocal part that occurs since we are considering non-singular kernels. 
\end{abstract}

\maketitle


\section{Introduction}

Nonlocal models can be used to describe phenomena that involve
long-range interactions and discontinuities. 
For example, long-range interactions effectively
describe anomalous diffusion and crack
formation in material models. 
The fundamental difference between nonlocal and classical local models  
is the fact that the latter only involve differential operators (local equations), whereas the
former rely on integral operators (nonlocal equations). For general references on nonlocal models 
and their applications, we refer to
to \cite{BCh,Bere,CF,ChChRo,CERW,Cortazar-Elgueta-Rossi-Wolanski,F,Hutson} and the book \cite{ElLibro}.

It is often the case that nonlocal effects occur only in
some parts of the domain, whereas, in the remaining parts, the physical setting can be accurately
described by a local equation. The goal of coupling local and nonlocal models is to
combine a local equation (a partial differential equation) with a nonlocal one (an integral equation) acting inside different
parts of the domain. 
In this
context, one of the challenges of a coupling strategy is to provide a mathematically 
consistent formulation.  Hence, interesting problems arise from coupling local and nonlocal models, 
see \cite{delia2,delia3,Gal,ROSSICVPDE,Han,Kri} and references therein. Previous strategies treat the coupling condition as an optimization
problem (the
goal is to minimize the mismatch of the local and nonlocal solutions in a common overlapping
region). Another possible strategy for coupling relies on the partitioned procedure as a general coupling
strategy for heterogeneous systems, the system is divided into sub-problems in
their respective sub-domains, which communicate with each other via transmission
conditions. Moreover, couplings between sets of different dimension are possible.
In \cite{Bere} the effects of network transportation on enhancing biological invasion is studied. 
The proposed mathematical model consists of one equation with nonlocal diffusion in a one-dimensional domain coupled via the boundary condition with a standard reaction-diffusion in a two-dimensional domain. 
In \cite{delia2}, local and nonlocal problems were coupled through a prescribed region in which both kinds of equations overlap (the value of the solution in the nonlocal part of the domain is used as a Dirichlet boundary condition for the local part and vice-versa). This kind of coupling gives continuity of the solution in the overlapping region but does not preserve the total mass when Neumann boundary conditions are imposed.
In \cite{Gal,ROSSICVPDE,Kri}, 
evolution problems related to energies closely related to ours are studied (here we deal with stationary problems).

Here we continue the analysis of coupled local/nonlocal models and deal with the validity of the maximum 
principle and the existence of a principal eigenvalue for these operators. 

\subsection*{Problem and Hypotheses}
 
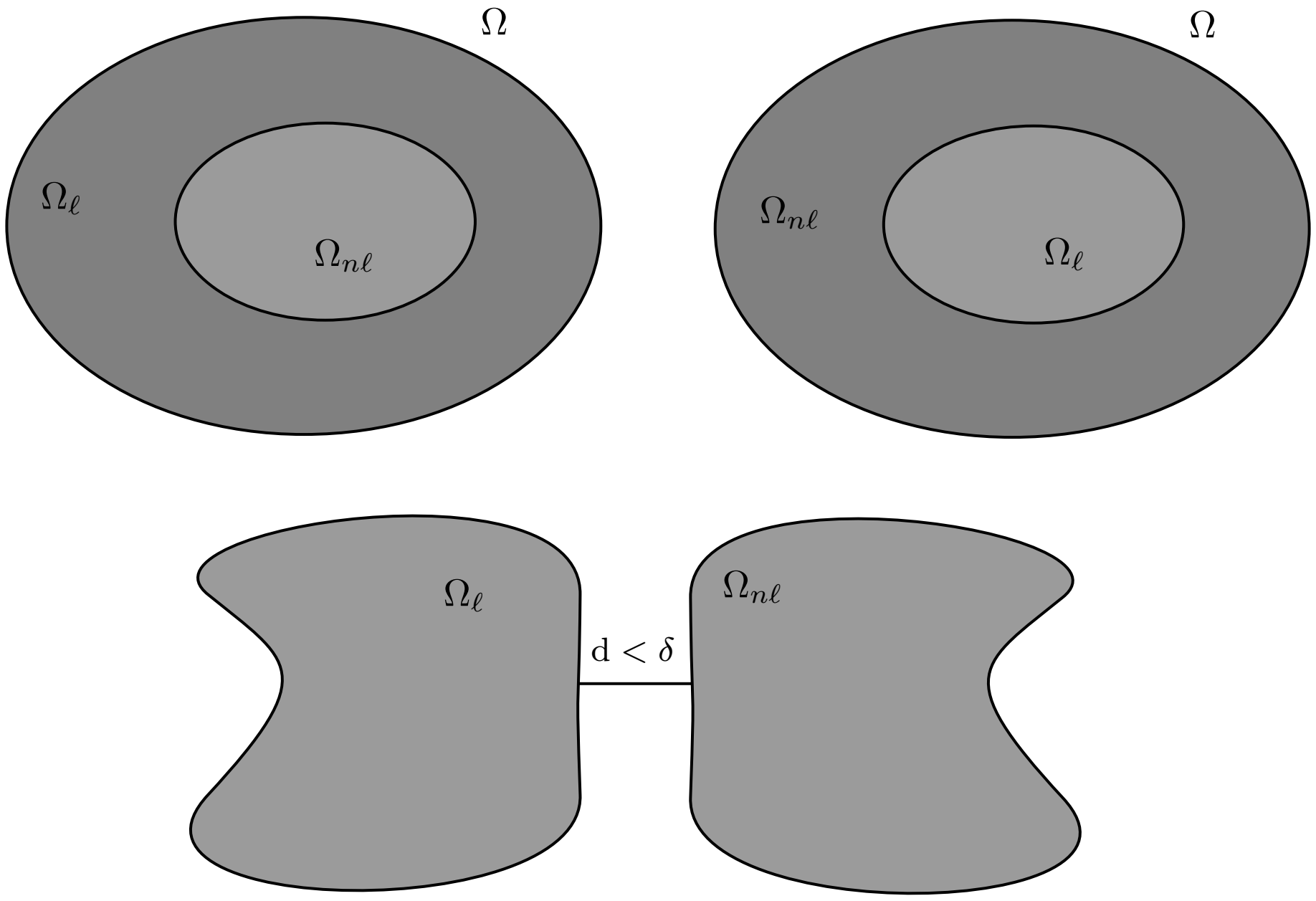
\begin{figure}
\centering\input{fig1}\caption{Admissible configurations of the domain. In the last case  $d:=\operatorname{dist}(\Oe,\Onl)<\delta$.}\label{fig: Dominio}\end{figure}

Let us describe more precisely the problem under consideration. 
Let $\Omega_{\ell}$ and $\Omega_{n\ell}$ be disjoint regular subdomains of $\mathbb{R}^N$ ($N \geq 1$). We define the total domain as $\Omega = \Omega_{\ell} \cup \Omega_{n\ell} \cup \Gamma$, where $\Gamma := \partial \Omega_{\ell} \cap \partial \Omega_{n\ell}$ denotes their common boundary. We assume that $\operatorname{dist}(\Omega_{\ell}, \Omega_{n\ell}) < \delta$ for a given constant $\delta > 0$. Note that if the subdomains are strictly separated, $\Gamma$ is empty. However, this formulation naturally incorporates the scenario where $\Omega$ is partitioned by a hypersurface $\Gamma$; in this case, $\operatorname{dist}(\Omega_{\ell}, \Omega_{n\ell}) = 0$, which trivially satisfies the threshold condition (see Figure \ref{fig: Dominio}).

In this paper, motivated by the model introduced by \cite{RossiSIAM, ROSSICVPDE} (see also  \cite{CAPANNAROSSIMJM, dasilva2026parabolicellipticdynamicslocalnonlocalcoupled}), we study the following eigenvalue problem
\[ \begin{cases}\tag{$P_{\lambda}$}
\Lop_1 u (x) + B_1(x)u (x) = \lambda u (x)  &x\in \Omega_{\ell}, \\  \mathfrak{B}u (x) = 0 & x\in \partial \Omega_{\ell}, \\  
\Lop_2 u (x) + B_2(x)u (x) =\lambda u  (x)\qquad & x \in  \Omega_{n \ell}, \label{P_l} 
\end{cases}
\] where $B_1 \in L^\infty(\Omega_\ell), B_2\in C(\overline{\Omega}_{n\ell})$. To rigorously define the operators involved, let us introduce the Banach space 
$$X= \{u \in L^2(\Omega) : u\vert{}_{\Omega_\ell} \in W^{2,p}(\Omega_\ell),\, p > N, \text{ and }  u\vert{}_{\overline{\Omega}_{n\ell}} \in C(\overline{\Omega}_{n\ell}) \}.$$
Thus, $\mathcal{L}_1: X \to L^2(\Omega_\ell)$ and $\mathcal{L}_2: X \to L^2(\Omega_{n\ell})$ are given by
\begin{align*}
\mathcal{L}_1u(x) &:= -\Delta u (x) - \int_{\Omega_{n\ell}} J(x-y)(u(y)-u(x))dy \\ 
\mathcal{L}_2u(x) &:= -2\int_{\Omega_{n\ell}} J(x-y)(u(y)-u(x))dy - \int_{\Oe} J(x-y)(u(y)-u(x))dy.
\end{align*}
Furthermore, for the dispersion kernel $J$ we assume that 
\begin{enumerate}[label=$(J_{\arabic*})$]\item \label{J1} $J \in C(\mathbb{R}^N)$ is a nonnegative function such that $J(0)>0$ and there exists $\delta>0$ such that $\operatorname{supp}J = B_\delta(0)$.
\end{enumerate}

Remark that in $\Lop_1$ that is the main part of the equation in $\Omega_{\ell}$ the main operator 
that appears in $-\Delta$ the most studied local operator; while in $\Lop_2 $ that acts in $\Omega_{n \ell}$ 
the main part is given by $-2\int_{\Omega_{n\ell}} J(x-y)(u(y)-u(x))dy $
that is a nonlocal operator. The equations that occur in these two parts of the domain are 
coupled by nonlocal transmission terms, $- \int_{\Omega_{n\ell}} J(x-y)(u(y)-u(x))dy$ in $\Omega_{\ell}$
and $ - \int_{\Oe} J(x-y)(u(y)-u(x))dy$ in $\Omega_{n \ell}$, that encode the total amount of mass that is
transferred from one subdomain to the other.  Also notice that in our eigenvalue problem we have two functions $B_1$ and $B_2$
and for them we are not assuming a priori any sign condition.

%The right hand side of the equations in \eqref{P_l} has a natural energy functional
%$$
%\begin{array}{l}
%\displaystyle 
%F(u) = \frac12 \int_{\Omega_{\ell}} |\nabla u(x) |^2 dx +  \frac12\int_{\Omega_{n \ell}} \int_\Omega J(x-y) (u(x)- u(y))^2 dx dy
%\\[6pt]
%\qquad\qquad \displaystyle + \frac12 \int_{\Omega_{\ell}} B_1 (x) (u(x))^2 dx + + \frac12 \int_{\Omega_{n\ell}} B_2 (x) (u(x))^2 dx.
%\end{array}
%$$

Under the assumptions of this work, $\partial \Omega_\ell \neq \emptyset$ is composed of two disjoint components, $\Gamma_{0}:=\partial \Omega_{\ell} \cap \partial\Omega$ and $\Gamma=\partial \Omega_{\ell} \cap \Omega$, either of which may be empty. Therefore, the boundary operator $\mathfrak{B} : C(\Gamma_0) \times C^1(\Gamma) \to C(\partial \Omega_\ell)$ is defined by
\begin{equation*}
\mathfrak{B} u := \begin{cases}u & \text{on } \Gamma_0, \\ \frac{\partial u}{\partial \eta} & \text{on } \Gamma,\end{cases}
\end{equation*} where $\eta=\eta(x)$ denotes the outward unit normal field to $\Gamma=\partial\Oe \cap \Omega$.  

We also introduce the following Banach spaces:
$$X_{\mathfrak{B}}= \{u \in X : \, \mathfrak{B}u=0 \text{ on } \partial \Omega_{\ell} \},$$
$$W_\mathfrak{B}^{2,p}(\Omega_{\ell}) = \{ u \in W^{2,p}(\Omega_{\ell}); \mathfrak{B}u=0\text{ on } \partial \Omega_{\ell}\},
$$ 
and 
$$Y=\{u \in L^2(\Omega); \, u\vert{}_{\Oe} \in L^p(\Oe), p>N, \text{ and } u\vert{}_{\Onl} \in C(\overline{\Omega}_{n\ell})\}.$$
\begin{remark}{\rm Note that we can make the natural identification of $X_\mathfrak{B}$ with $W_\mathfrak{B}^{2,p}(\Omega_\ell) \times C(\overline{\Omega}_{n\ell})$ such that each $w \in X_\mathfrak{B}$ is identified with the pair $(u,v)\in W_\mathfrak{B}^{2,p}(\Omega_\ell) \times C(\overline{\Omega}_{n\ell})$ where $u=w\vert{}_{\Oe}$ and $v=w\vert{}_{\Onl}$. Analogously, $Y$ is identified with $L^p(\Oe)\times C(\overline{\Omega}_{n\ell})$.}
\end{remark}

With these notations, we define $B:=(B_1,B_2)$ and the operator $L:X \to Y$ given by
$$Lu:= \begin{cases} \Lop_1 u  & \text{in } \Omega_{\ell}, \\ \Lop_2 u  & \text{in } \Omega_{n \ell}. \end{cases}$$
Thus, the eigenvalue problem \eqref{P_l}  can be rewritten as
$$ \begin{cases}  Lu (x) +B(x)u (x) = \lambda u (x) & x\in  \Omega,\\  
\mathfrak{B} u (x) =  0& x \in \partial \Omega_{\ell}.  \end{cases}$$

We will look for solutions of \eqref{P_l}  in the strong sense. That is, a function $u \in X_\mathfrak{B}$ such that $u\vert{}_{\Omega_\ell}\in W^{2,p}_\mathfrak{B}(\Omega_\ell)$ satisfies the differential equation in $\Omega_\ell$ almost everywhere, while $u\vert{}_{\Omega_{n\ell}}\in C(\overline{\Omega}_{n\ell})$ satisfies the integral equation pointwise.

Note that if $B\equiv 0$ and $\Gamma_0 = \emptyset$, then $u=1$ is a positive solution of \eqref{P_l} for $\lambda=0$ (that is, $\lambda=0$ is a principal eigenvalue). However, under the general assumptions considered in this paper, \eqref{P_l} might not admit a positive solution for any $\lambda \in \mathbb{R}$, as happens in the purely nonlocal case (see, for instance, \cite[Section 5]{MR2718672}).

\begin{remark}
  {\rm  We point out that under condition \ref{J1}, the integral coupling terms in \eqref{P_l} vanish identically whenever $\operatorname{dist}(\Omega_{\ell}, \Omega_{n\ell}) \geq \delta$, thus decoupling \eqref{P_l} into the following two classical eigenvalue problems
    $$
    \left\{
    \begin{array}{ll}
    -\Delta u (x)+B_1(x)u (x) = \lambda u (x)   & x\in \Oe,  \\
    \mathfrak{B} u (x) =0     & x\in \partial \Oe,
    \end{array} \right.
    $$
    and 
    $$  -2 \int_{\Onl} J(x-y) (u(y)-u(x)) \, dy +B_2(x)u (x) = \lambda u (x) \quad \mbox{in }\Onl.$$
    
    For these problems we refer to \cite{BNV,bookSMP} and to \cite{BCV,MR2718672,Jorge,MR3583503}. 
}
\end{remark}

It is worth recalling that the analysis of the principal eigenvalue problem and its properties is fundamental to the understanding of nonlinear problems using, for instance, sub- and supersolution methods and bifurcation theorems. Moreover, the validity of the maximum principle is crucial when analyzing positive solutions. For this reason, the current work is the initial step towards investigating coupled local/nonlocal nonlinear stationary problems. Additionally, the lack of compactness of the nonlocal operator makes its mathematical treatment more challenging. For instance, the approach via Krein-Rutman theorems adopted in \cite{interfaceJMAA} for purely local coupled problems is not applicable. To overcome this difficulty, inspired by \cite{BCV,BNV,MR2718672,DV75,MR3583503}, we employ a pointwise min-max characterization associated with \eqref{P_l} and gather several properties by combining results involving both local and nonlocal operators.

\subsection*{Notations and Main results:} Let us 
introduce some notations that will be used in the sequel,
\begin{itemize}\item $P_\ell:=P_{C_{0}^{1}(\overline{\Omega}_{\ell})}$ and $P_{n\ell}:=P_{C(\overline{\Omega}_{n\ell})}$ are the positive cones of $C_{0}^{1}(\overline{\Omega}_{\ell})$ and $C(\overline{\Omega}_{n\ell})$ respectively.
\item When $u \in X$, we write $u\gg 0$  to mean that
$$u\vert{}_{\Omega_\ell} \in \operatorname{int} P_{\ell} \quad \mbox{and} \quad u\vert{}_{\overline{\Omega}_{n\ell}} \in \operatorname{int} P_{n\ell}.$$
\item When $u \in Y$, we write $u>0$  to mean that $u \geq 0$ a.e. in $\Omega$ and $u \neq 0$. In this case we say that $u$ is {\it positive}. Similarly for $u>0$ in $\Onl$ or $u>0$ in $\Oe$.
\item $\mathfrak{J}_{n\ell}(x):=\int_{\Omega_{n \ell}}J(x-y)dy$ and $\mathfrak{J}_{\ell}(x):=\int_{\Oe}J(x-y)dy.$\item $\vert{}\cdot\vert{}_{2}$ denotes the usual norm in $L^{2}(\Omega)$.\item $\langle\cdot,\cdot\rangle_2$ denotes the usual inner product in $L^2(\Omega)$.
\item $\rho(L+B)$ denotes the resolvent of the operator $L+B$.
\item $\sigma(L+B):=\mathbb{R}\setminus \rho(L+B)$ denotes the spectrum of the operator $L+B$.
\end{itemize}

The paper is organized as follows, outlining our main contributions by section:

Section \ref{sect-Maximum} shows that the existence of a positive strict supersolution of operator $(L+B,\mathfrak{B})$ in $\Omega$ guarantees the validity of the strong maximum principle (see Corollary \ref{SS-SMP}), and we also deduce some important consequences of this result.

In Section \ref{sect-Eigen}, we introduce the set
$$ \mathfrak{P}= \{ u \in X:  u(x)>0 \mbox{ for all }x \in \Omega \mbox{ and } \mathfrak{B}u \geq 0 \mbox{ on }\partial\Oe\}$$
and define
$$\lambda_{p} := \sup_{\phi\in\mathfrak{P}} \inf_{x \in \Omega} \frac{(L+B(x))\phi(x)}{\phi(x)}.$$
First, it follows that
$$ \lambda_p \leq m := \min_{x \in {\overline{\Omega}_{n\ell}}} \left(
2 \int_{\Onl} J(x-y)dy + \int_{\Oe} J(x-y)dy + B_2(x) \right).$$
We then prove that $\lambda_p$ is a principal eigenvalue of \eqref{P_l} (in the sense that it admits a positive eigenfunction) if and only if $\lambda_p < m$. Furthermore, whenever this condition holds, we provide several key properties of this eigenvalue and its corresponding eigenfunction (see Proposition \ref{Prop:lambda-av-principal} and Theorem \ref{Th: caso l_p<m}). We also present a sufficient condition to ensure that $\lambda_p < m$.

Section \ref{sect-proper} is devoted to formulating useful equivalent definitions for $\lambda_p$, as well as deriving various monotonicity properties of this value. Specifically, we show that
$$
\lambda_p=\inf_{\substack{\phi \in X_{\mathfrak{B}}\\ \phi>0}} \sup_{x \in \Omega} \frac{(L+B(x))\phi(x)}{\phi(x)}=\sup \mathcal{A}=\inf \Tilde{\mathcal{A}},
$$
where
$$
\mathcal{A} := \Big\{ \lambda \in \mathbb{R} : \exists \, \phi \in X, \, \phi(x) > 0, \, x \in \Omega,
     \text{ supersolution of } (L + B - \lambda, \mathfrak{B}, \Omega) \Big\},
$$
 and
$$
   \Tilde{\mathcal{A}} := \Big\{ \lambda : \exists \, \phi \in X_{\mathfrak{B}}, \, \phi > 0, \, \text{subsolution of } (L+ B-\lambda, \mathfrak{B}, \Omega) \Big\}.
$$

Finally, Section \ref{sect-variational} assumes the symmetry of $J$ to obtain the following variational characterization for $\lambda_p$,
$$
    \lambda_{p} =\inf_{\substack{\phi \in H\\ |\phi|_2=1}} \displaystyle \left(\int_{\Oe} |\nabla \phi |^{2}dx+\int_{\Onl}\int_{\Omega}J(x-y)(\phi(y)-\phi(x))^{2}dy dx+\int_{\Omega}B(x)\phi^2(x)dx\right)
  $$
  with 
  $$H= \Big\{ u \in L^2(\Omega); u|_{\Oe} \in H^1_{\Gamma_0}(\Oe)\Big\}.$$

\section{Maximum Principle Results} \label{sect-Maximum}

This section is devoted to the maximum principle for the operator $(L+B, \mathfrak{B}, \Omega)$. We introduce the notion of supersolutions, provide a classification result, and establish sufficient conditions under which the strong maximum principle holds. Furthermore, we discuss some important consequences of these results. 

We begin with the following fundamental definitions:
\begin{definition}
We say that $u \in X$ is a supersolution of $(L+B, \mathfrak{B}, \Omega)$ if
$$\begin{cases} 
Lu (x)+B(x)u (x) \ge 0 & x\in \Omega, \\ 
\mathfrak{B} u (x) \ge 0 & x\in  \partial \Omega_\ell.
\end{cases}$$
When one of the above inequalities is strict, we say that $u$ is a strict supersolution.
\end{definition}

 Similarly, we define a subsolution of $(L+B, \mathfrak{B}, \Omega)$ by reversing the inequalities.

\begin{definition}
We say that $(L+B, \mathfrak{B}, \Omega)$ satisfies the Strong Maximum Principle if every supersolution $u\in  X \setminus \{0\}$ of $(L+B, \mathfrak{B}, \Omega)$ satisfies $u \gg0$.
\end{definition}

The following result establishes a counterpart to the second classification of supersolutions from \cite[Theorem 2.4]{bookSMP} or the Touching Lemma in \cite[Lemma 2.1]{MR3583503}, for the problem considered in this work.

\begin{theorem}\label{Teorema C}Suppose that $(L+B, \mathfrak{B}, \Omega)$ admits a positive supersolution $h \in X$. Then, every supersolution $u \in X$ satisfies one of the following alternatives:
\begin{enumerate}[label=$(A_{\arabic*})$]
    \item \label{A_1}$u \equiv0$;
    \item \label{A_2} $u \gg 0$;
    \item \label{A_3}There exists a constant $\tau < 0$ such that $u = \tau h$. In this case,
    \[\begin{cases}
        Lu (x)+B(x)u (x)=0 \ & x\in   \Omega, \\
        \mathfrak{B}u (x)=0 & x\in   \partial\Oe.
    \end{cases}\]
\end{enumerate}
\end{theorem}

\begin{proof} First, we will show that $h \gg 0$. Since $h$ is positive, it follows that either $h>0$ in $\Oe$ or $h>0$ in $\Onl$. We analyze these two possibilities:

\textit{\textbf{Case 1}:} Suppose that $h>0$ in $\Oe$ and let $M>0$ be a constant such that $\mathfrak{J}_{n\ell} +B_1+ M >0$ in $\Oe$. Since $h>0$ is a supersolution of $(L+B,\mathfrak{B},\Omega)$, we have
$$\begin{cases}
    -\Delta h+h(\mathfrak{J}_{n\ell}(x)+B_1(x)+M)\ge \Lop_1 h + (B_1(x)+M)h \ge Mh> 0 & \text{in } \Omega_{\ell}, \\
\mathfrak{B}h \ge 0 &  \text{on } \partial \Oe. 
\end{cases}$$
By the choice of $M$, the local operator $(-\Delta +\mathfrak{J}_{n\ell}+B_1+M, \mathfrak{B}, \Oe)$ satisfies the Strong Maximum Principle (see, for instance, \cite{bookSMP}). Consequently, $h|_{\Omega_\ell} \in \operatorname{int}P_\ell$. Next, we show that $h \not\equiv 0$ in $\Onl$. Indeed, suppose by contradiction that $h \equiv 0$ in $\Onl$. Using again the fact that $h$ is a supersolution of $(L+B,\mathfrak{B},\Omega)$, we conclude that
\begin{align} \label{2.1}\mathcal{L}_2h = -\int_{\Oe} J(x-y)h(y)dy \ge 0 \quad \text{in } \Onl.\end{align}
On the other hand, as long as $\operatorname{dist}(\Omega_{\ell}, \Omega_{n\ell}) < \delta$, we can find $x_0 \in \Onl$ such that
$$A := B_{\delta}(x_0) \cap \Oe \ne \emptyset.$$
By hypothesis \ref{J1}, and since $h \in \operatorname{int}P_\ell$, it follows that
$$\int_{\Oe} J(x_0-y)h(y)dy = \int_{A} J(x_0-y)h(y)dy > 0,$$
which contradicts \eqref{2.1}. This shows that $h|_{\Onl} > 0$. Thus, we have
$$H_+ := \{x \in \Omega_{n\ell} : h(x) > 0\} \ne \emptyset.$$
We claim that $H_+ = \Omega_{n\ell}$. Otherwise, for $x_1 \in \partial H_+ \subset \Onl$, we would have
\begin{align*}0 \leq \mathcal{L}_{2}h(x_1) + B_2(x_1)h(x_1) = -2\int_{\Onl} J(x_1-y)h(y)dy - \int_{\Oe} J(x_1-y)h(y)dy < 0,\end{align*}
which is impossible.  Moreover, if $h(x^*)=0$ for some $x^*\in\partial \Onl$, by choosing $x_k \in \Onl$ such that $x_k \to x^*$, then
\begin{multline*}
    0 \leq \Lop_2 h(x_k) + B_2(x_k)h(x_k)  \to \Lop_2 h(x^*) + B_2(x^*)h(x^*)\\ = -2\int_{\Onl} J(x^*-y)h(y)dy - \int_{\Oe} J(x^*-y)h(y)dy < 0.
\end{multline*}
Consequently, $h(x)>0$ for all $x \in \overline{\Omega}_{n\ell}$, that is, $h|_{\Onl} \in \operatorname{int}P_{n\ell}$.

\textit{\textbf{Case 2}:} Now suppose that $h > 0$ in $\Onl$. 
In this case, we again have $H_+ := \{x \in {\Omega}_{n\ell} : h(x) > 0\} \ne \emptyset$, and by the same arguments as above, we conclude that $h|_{\overline{\Omega}_{n\ell}} \in \operatorname{int}P_{n\ell}$. In particular, it follows that 
\[
 \int_{\Onl} J(x-y) h(y) dy > 0,
\]
for all $x \in \Oe$ such that $\operatorname{dist}(x,\Onl)<\delta$. Therefore, since $h$ is a supersolution, we obtain
\[
\begin{cases}
   \displaystyle  -\Delta h + (\mathfrak{J}_{n\ell}(x) + B_{1}(x))h \ge \int_{\Onl} J(x-y) h(y) dy > 0 & \text{in } \Oe, \\
    \mathfrak{B} h \ge 0 & \text{on } \partial \Oe.
\end{cases}
\]
That is, $h|_{\Oe}$ is a positive supersolution for the local operator $(-\Delta + \mathfrak{J}_{n\ell} + B_{1}, \mathfrak{B}, \Oe)$. By classical maximum principle results, it follows that $h|_{\Oe} \in \operatorname{int} P_\ell$. Thus, we ensure that $h \gg 0$.

Now, let $u \in X$ be a supersolution. If $u \equiv 0$, there is nothing to prove. If $u > 0$ in $\Omega$, we can repeat the arguments used for the function $h$ and conclude that $u$ satisfies \ref{A_2}. It remains to consider the case where there exists $x_- \in \Omega \cup \Gamma$ such that $u(x_-) < 0$. Let us define the set:
\[
\Lambda := \{t > 0 : v_{t} := u + th \ge 0 \text{ in } \Omega \}.
\]
Since $t \in \Lambda$ for sufficiently large $t$, then $\mu = \inf \Lambda$ is well-defined. Furthermore, it is clear that $\mu > 0$ because $u(x_-) < 0$. Thus, we have $v_\mu \geq 0$. If $v_\mu \equiv 0$, then \ref{A_3} holds with $\tau = -\mu$. On the other hand, if $v_\mu > 0$, we can once again repeat the arguments made for the function $h$ and conclude that $v_\mu \gg 0$. However, this implies the existence of some $\mu' < \mu$ such that $v_{\mu'} \geq 0$, contradicting the minimality of $\mu$. This completes the proof.
\end{proof}

The preceding theorem allows us to obtain a sufficient condition for the operator $(L+B, \mathfrak{B}, \Omega)$ to satisfy the Strong Maximum Principle.

\begin{corollary}
\label{SS-SMP}
If $(L+B, \mathfrak{B}, \Omega)$ admits a positive strict supersolution $h \in X$, then $(L+B, \mathfrak{B}, \Omega)$ satisfies the Strong Maximum Principle.
\end{corollary}

\begin{proof}
Suppose that $h \in X$ is a positive strict supersolution of $(L+B, \mathfrak{B}, \Omega)$. Arguing as in the proof of Theorem \ref{Teorema C}, we have $h \gg 0$. Now, let $u \in X \setminus \{0\}$ be a supersolution of $(L+B, \mathfrak{B}, \Omega)$. Again by Theorem \ref{Teorema C}, we have that $u$ satisfies one of the alternatives \ref{A_1}--\ref{A_3}. Since $u \neq 0$, \ref{A_1} cannot occur. To conclude, it remains to show that \ref{A_3} also cannot be satisfied. Indeed, suppose by contradiction that there exists a constant $\tau < 0$ such that $u = \tau h$. Thus,
$$
\left\{\begin{array}{ll}
     0 \leq Lu+B u = \tau (Lh+B(x)h) \leq 0 &  \text{in } \Omega, \\
     0 \leq \mathfrak{B}u = \tau \mathfrak{B}h \leq 0 & \text{on } \partial \Oe. 
\end{array} \right.
$$
Since $h$ is a strict supersolution, one of the inequalities above must be strict, which is impossible. Consequently, $u$ satisfies \ref{A_2} and the result follows.
\end{proof}

We conclude this section with two consequences of Corollary \ref{SS-SMP}.

\begin{corollary} \label{Cor-nucleo-trivial}
Suppose that $(L+B, \mathfrak{B}, \Omega)$ admits a positive strict supersolution. Then the homogeneous linear problem
\begin{equation}\label{P0}
    \begin{cases}
    L u (x) + B(x)u (x) = 0 &x \in\Omega, \\
\mathfrak{ B}u (x) = 0 &x \in \partial \Oe, \\
\end{cases}
\end{equation}
possesses only the trivial solution $u\equiv0$. That is, $\operatorname{ker}[L+B] = \{0\}$, where $L+B: X_\mathfrak{B} \to Y$.
\end{corollary}

\begin{proof}
Suppose, by contradiction, that there exists $u \in X \setminus \{0\}$ satisfying \eqref{P0}. In this case, $u \neq 0$ is a supersolution of $(L+B, \mathfrak{B}, \Omega)$. By Corollary \ref{SS-SMP}, it follows that $u \gg 0$. On the other hand, $-u \neq 0$ is also a supersolution of $(L+B, \mathfrak{B}, \Omega)$; thus, Corollary \ref{SS-SMP} again implies that $-u \gg 0$, which is a contradiction. This completes the proof.
\end{proof}

\begin{corollary}
Suppose that $(L+B, \mathfrak{B}, \Omega)$ admits a positive strict supersolution. Then, given $f \in L^2(\Omega)$ and $g \in L^2(\partial \Omega_\ell)$, the linear problem
\begin{equation}\label{PL}
    \begin{cases}
L u (x) + B(x)u (x) = f(x) & x\in \Omega, \\
\mathfrak{B} u (x) = g(x) & x\in \partial \Oe, \\
\end{cases}
\end{equation}
admits at most one solution in $X$.
\end{corollary}

\begin{proof}
It suffices to observe that if $u, v \in X$ are solutions of \eqref{PL}, then $u-v \in X_\mathfrak{B}$ is a solution of \eqref{P0}, and by the previous Corollary, we have $u-v=0$.
\end{proof}

\section{Eigenvalue Problem} \label{sect-Eigen}

This section is devoted to the study of the eigenvalue problem \eqref{P_l}. Our goal is to provide sufficient conditions for this problem to admit a positive solution, as well as to study several qualitative properties. We begin with the following definition.

\begin{definition}
    We say that a value $\lambda \in \mathbb{R}$ is an eigenvalue of the operator $(L+B,\mathfrak{B},\Omega)$ if the problem 
    \[
\begin{cases}
    L u (x)+B(x)u (x) = \lambda u (x) & x\in \Omega, \\
\mathfrak{ B}u (x) = 0 & x\in \partial \Omega_\ell, \\
\end{cases}
\]
admits a nontrivial solution in $X_\mathfrak{B}$. Additionally, we say that $\lambda \in \mathbb{R}$ is a principal eigenvalue if it admits a nontrivial solution in $X_\mathfrak{B}$ that does not change sign (notice that we can always assume that the eigenfunction associated with a principal eigenvalue is positive).
\end{definition}

Based on the previous definition, and motivated by the works of \cite{BNV,MR2718672,DV75}, we introduce the following definitions:
$$
\mathfrak{P}= \Big\{ u \in X:  u(x)>0 \mbox{ for all }x \in \Omega \mbox{ and } \mathfrak{B}u \geq 0 \mbox{ on }\partial\Oe \Big\}
$$
and
\begin{equation}
\lambda_{p} := \sup_{\phi\in\mathfrak{P}} \inf_{x \in \Omega} \frac{(L+B(x))\phi(x)}{\phi(x)}. \label{Def lambdap}
\end{equation}

\begin{remark}\label{Remark-l<lp} {\rm
    Note that, by definition, for every $\lambda < \lambda_p$, there exists $\phi_\lambda \in \mathfrak{P}$  such that
    $$
    \lambda < \inf_{\Omega} \frac{(L+B(x))\phi_\lambda(x)}{\phi_\lambda(x)},
    $$
    which implies that 
\begin{equation*}
    \begin{cases}
    L \phi_\lambda (x)+ B(x)\phi_\lambda (x) - \lambda \phi_\lambda (x) > 0 & x \in \Omega, \\
\mathfrak{ B}\phi_\lambda (x) \ge 0 & x \in \partial \Oe. \\
\end{cases}
\end{equation*}
That is, $\phi_\lambda > 0$ is a positive strict supersolution for $(L+B-\lambda,\mathfrak{B},\Omega)$. In particular, this operator satisfies the Strong Maximum Principle (see Corollary \ref{SS-SMP}) and $\operatorname{ker}[L+B-\lambda]=\{0\}$ (see Corollary \ref{Cor-nucleo-trivial}). }
\end{remark}

From now on, we use the following notation:
\[
a(x) := 2 \int_{\Onl} J(x-y)dy + \int_{\Oe} J(x-y)dy + B_2(x) 
\]
and
$$
m := \min_{x \in {\overline{\Omega}_{n\ell}}} a(x).
$$
It is also convenient to introduce
$$
L_0u(x):= -2\int_{\Onl} J(x-y)u(y)dy - \int_{\Oe} J(x-y)u(y)dy \quad x \in \Onl.
$$
Thus, by a direct calculation, it can be verified that  
\begin{align}
\Lop_2 u (x) + B_2(x)u (x) &= L_0u (x) + a(x)u (x) \quad \mbox{in }\Omega_{n\ell}. \label{Equ. nl}
\end{align}

\begin{remark}\label{Remark lambdap menor m} {\rm 
     Using \eqref{Equ. nl}, for every $\phi \in \mathfrak{P}$, we have
$$\frac{(L+B(x))\phi(x)}{\phi(x)}=\frac{(\Lop_2  + B_2(x))\phi(x)}{\phi(x)}  \leq  a(x) \quad \mbox{in }\Onl. $$
Thus, it follows from the definition  \eqref{Def lambdap} that $$\lambda_p \leq m. $$}
\end{remark}

With the above considerations, we provide a necessary condition for $\lambda \in \mathbb{R}$ to be a principal eigenvalue of $(L+B,\mathfrak{B},\Omega)$.

\begin{proposition}\label{Prop:lambda-av-principal}
If $\lambda \in \mathbb{R}$ is a principal eigenvalue of $(L+B,\mathfrak{B},\Omega)$, then $\lambda_p \leq \lambda < m$.
\end{proposition}

\begin{proof}
Let $\varphi \in X_\mathfrak{B}$ with $\varphi > 0$ be an eigenfunction associated with $\lambda$. Then $\varphi$ is a positive supersolution of $(L+B-\lambda,\mathfrak{B},\Omega)$, and it follows from the proof of Theorem \ref{SS-SMP} that $\varphi \gg 0$. Let  $x_0 \in \overline{\Omega}_{n\ell}$ with $a(x_0) = m$ and $x_n \in \Omega_{n\ell}$ a sequence such that $x_n \to x_0$. Using \eqref{Equ. nl} we obtain 
\begin{multline*}
    \lambda \varphi(x_n) = \mathcal{L}_2\varphi(x_n) + B_2(x_n)\varphi(x_n) = \cr -2\int_{\Onl} J(x_n-y)\varphi(y)dy - \int_{\Oe} J(x_n-y)\varphi(y)dy + a(x_n)\varphi(x_n).
\end{multline*}
Letting $n \to +\infty$ we find that
$$
    \lambda \varphi(x_0) = -2\int_{\Onl} J(x_0-y)\varphi(y)dy - \int_{\Oe} J(x_0-y)\varphi(y)dy + a(x_0)\varphi(x_0)< m\varphi(x_0).
$$
Since $\varphi(x_0) > 0$ (because $\varphi|_{\overline{\Omega}_{n\ell}} \in \operatorname{int}P_{n\ell}$), we conclude that $\lambda < m$.

On the other hand, if $\lambda < \lambda_p$, it follows from Remark \ref{Remark-l<lp} that $\operatorname{ker}[L+B-\lambda]=\{0\}$. Therefore, $\lambda$ cannot be an eigenvalue of $(L+B,\mathfrak{B},\Omega)$.
\end{proof}

The following result shows that the condition $\lambda_p<m$ is sufficient for $\lambda_p$ to be a principal eigenvalue. The proof is inspired by \cite{MR3583503}.

\begin{theorem}\label{Th: caso l_p<m}
Assume that $\lambda_p < m$. Then the following statements are true.
\begin{enumerate}
\item[(i)] $\lambda_p$ is a principal eigenvalue of $(L+B,\mathfrak{B},\Omega)$ that admits an eigenfunction $\varphi_p \in X_\mathfrak{B}$ such that $\varphi_p \gg0$.
\item[(ii)] $\lambda_p$ is the unique principal eigenvalue of $(L+B,\mathfrak{B},\Omega)$.
\item[(iii)] If $J$ is symmetric then $\lambda_{p}$ is isolated in the spectrum of $\sigma(L+B)$ and  simple, in the sense that $\operatorname{ker}[L+B-\lambda_pI]=\operatorname{span}\{\varphi_p\}$ and the equation
\begin{equation}\label{simple-eigenvalue}
    \begin{cases}
    L u (x) + B(x)u (x) - \lambda_p u (x) = \varphi_p (x)  & x\in \Omega, \\
\mathfrak{ B}u (x) = 0 & x\in \partial \Oe, \\
\end{cases}
\end{equation}
does not admit a solution $u \in X_\mathfrak{B}$. 
\end{enumerate}
\end{theorem}

\begin{proof}
We first prove that for every $\lambda < m$, the operator $L+B-\lambda: X_\mathfrak{B} \to Y$ is Fredholm of index zero. For this purpose, it is convenient to use the identifications $X_\mathfrak{B} = W_\mathfrak{B}^{2,p}(\Oe) \times C(\overline{\Omega}_{n\ell})$ and $Y = L^p(\Oe) \times C(\overline{\Omega}_{n\ell})$. Note that for each $\gamma \in \mathbb{R}$ and $u=(v,w) \in X_\mathfrak{B}$, we have
\begin{multline*}
     (L+B(x)-\lambda)(v,w) = \cr \begin{bmatrix}
        - \Delta v + (B_1(x)+ \mathfrak{J}_{n\ell}(x)+\gamma-\lambda)v - \displaystyle\int_{\Onl}J(x-y)w(y) \, dy - \gamma v \\
        -2\displaystyle\int_{\Onl} J(x-y)w(y)\,dy - \int_{\Oe} J(x-y)v(y)\,dy + (a(x)-\lambda)w
    \end{bmatrix}.
\end{multline*}

Thus, defining the operators $\Phi, K: X_\mathfrak{B} \to Y$ by
\[
\Phi = \begin{bmatrix}
- \Delta + (B_1(x) + \mathfrak{J}_{n\ell}(x) + \gamma - \lambda) & 0 \\[1.5em]
0 & (a(x)-\lambda)I
\end{bmatrix}
\]
and
\begin{align*}
K = \begin{bmatrix}
    -\gamma I & - \displaystyle\int_{\Onl}J(x-y) \cdot dy \\
    -\displaystyle\int_{\Oe} J(x-y) \cdot dy & -2\displaystyle\int_{\Onl} J(x-y) \cdot dy
\end{bmatrix},
\end{align*}
we have $L+B(x)-\lambda = \Phi + K$.

Next, we observe that by choosing $\gamma > 0$ sufficiently large, the principal eigenvalue of the (local) operator $(-\Delta + B_1 + \mathfrak{J}_{n\ell} + \gamma - \lambda, \mathfrak{B}, \Oe)$ is positive. This implies that $-\Delta + B_1 + \mathfrak{J}_{n\ell} + \gamma - \lambda: W_\mathfrak{B}^{2,p}(\Oe) \to L^p(\Oe)$ is a topological isomorphism. Furthermore, since $\lambda < m$, the multiplication operator $(a(x)-\lambda)I: C(\overline{\Omega}_{n\ell}) \to C(\overline{\Omega}_{n\ell})$ is also  invertible. Consequently, $\Phi: X_\mathfrak{B} \to Y$ constitutes an isomorphism between these spaces. 

On the other hand, from Sobolev compact embeddings, it follows that $\gamma I: W_\mathfrak{B}^{2,p}(\Oe) \to L^p(\Oe)$ is a compact operator. This fact, combined with the compactness of the integral operators, implies that $K$ is compact. Therefore, $L+B(x)-\lambda$ is the sum of an invertible operator and a compact one, and thus, it is a Fredholm operator of index zero.

Finally, recall that, as noted in Remark \ref{Remark-l<lp}, $\operatorname{ker}[L+B-\lambda] = \{0\}$ for all $\lambda < \lambda_p$. Since $L+B-\lambda: X_\mathfrak{B} \to Y$ is Fredholm of index zero, this trivial kernel implies that the operator is bijective, whenever $\lambda < \lambda_p$.

We now proceed to prove item (i). By the previous assertion, for any constant $\epsilon \in (0,1)$ and $\lambda \in [\lambda_p-1, \lambda_p)$, there exists a unique $\phi_\lambda^\epsilon \in X_{\mathfrak{B}}$ such that:
\begin{equation}\label{P_ep}
    \begin{cases}
(L+B-\lambda)\phi_\lambda^\epsilon(x) = \epsilon & x \in \Omega, \\
\mathfrak{B} \phi_\lambda^\epsilon(x) = 0 & x \in\partial \Oe. 
\end{cases}
\end{equation}
Since $(L+B-\lambda,\mathfrak{B},\Omega)$ satisfies the strong maximum principle (see Remark \ref{Remark-l<lp}), we have $\phi_\lambda^\epsilon(x) \gg 0$. We shall analyze the following cases:

\textbf{Case 1:} There exists $\epsilon \in (0,1)$ such that for all $\lambda \in [\lambda_p-1, \lambda_p)$, $|\phi_\lambda^\epsilon|_2 \le 1$.

\textbf{Case 2:} For every $\epsilon \in (0,1)$, there exists $\lambda_{\epsilon} \in [\lambda_{p}-1, \lambda_p)$ such that $|\phi_{\lambda_\epsilon}^{\epsilon}|_{2} > 1$.

We first show that \textbf{Case 1} leads to a contradiction. Indeed, from the definition of $L+B-\lambda$, it follows that
\begin{align}
\begin{cases}
\Lop_{1} \phi_{\lambda}^{\epsilon}+ B_1(x)\phi_\lambda^\epsilon(x) - \lambda \phi_\lambda^\epsilon(x) = \epsilon & x \in \Oe, \\
\mathfrak{B} \phi_\lambda^\epsilon(x) = 0 & x \in\partial \Oe. \label{2.26}
\end{cases}
\end{align}
Multiplying the first equation above by $\phi_\lambda^\epsilon$, integrating over $\Oe$, and using the fact that $|\phi_\lambda^\epsilon|_{2} \le 1$, we obtain
\[
\int_{\Oe} |\nabla \phi_\lambda^\epsilon|^2 \le \epsilon \int_{\Oe} \phi_\lambda^\epsilon - \int_{\Oe} (B_1(x)-\lambda)(\phi_\lambda^\epsilon)^2 + \int_{\Oe}\int_{ \Onl} J(x-y)\phi_\lambda^\epsilon(y)\phi_\lambda^\epsilon(x) dx dy < C,
\]
for some constant $C>0$. Hence, $\{\phi_\lambda^\epsilon\}$ is bounded in $H^1(\Oe)$. Thus, there exist a sequence $\lambda_n \uparrow \lambda_p$ and a function $\phi_\ell^\epsilon \in H^1(\Oe)$ such that
\begin{align*}
\phi_n^\epsilon := \phi_{\lambda_n}^\epsilon \rightharpoonup \phi_\ell^\epsilon \quad \text{in } H^1(\Oe), \\
\phi_n^\epsilon := \phi_{\lambda_n}^\epsilon \to \phi_\ell^\epsilon \quad \text{in } L^2(\Oe).
\end{align*}
Moreover, by the compactness of the integral operator and since $|\phi_n^{\epsilon}|_2 \leq 1$, we have
\[
\int_{\Onl} J(x-y)\phi_n^{\epsilon}(y)dy \to h \quad \text{in } C(\overline{\Omega}_{\ell}).
\]
Consequently, by a standard elliptic regularity argument, we conclude that
\[
\phi_n^\epsilon \to \phi_\ell^\epsilon \quad \text{in } W^{2,p}(\Oe), \text{ for } p>N.
\]
 In fact, note that 
\begin{align*}
    \begin{cases}
        -\Delta \phi_{n}^{\epsilon} = f_{n}^{\epsilon}(\lambda_n, \phi_{n}^{\epsilon}) \ &\text{in } \Oe, \\
        \mathfrak{B}\phi_{n}^{\epsilon} = 0 \ \ &\text{on } \partial \Oe,
    \end{cases}
\end{align*}
where
\[
f_{n}^{\epsilon}(\lambda_n,\phi_{n}^{\epsilon})
:= \epsilon+\lambda_n \phi_{n}^{\epsilon}
+ \int_{\Onl} J(x-y)\phi_{n}^{\epsilon}(y)\,dy
- \phi_{n}^{\epsilon}(x)\int_{\Onl} J(x-y)\,dy
- B_1(x)\phi_{n}^{\epsilon}.
\]
We observe that $f_{n}^{\epsilon}(\lambda_n, \phi_{n}^{\epsilon}) \in L^{2}(\Oe)$. Furthermore, the results established above show that this sequence converges strongly in $L^2(\Oe)$.
 Thus, by the Agmon-Douglis-Nirenberg Theorem, we have
\[
\|\phi_{n}^{\epsilon}-\phi_{m}^{\epsilon}\|_{W^{2,2}(\Oe)} 
\le C |f_{n}^{\epsilon}-f_{m}^{\epsilon}|_{L^{2}(\Oe)},
\]
which implies that $\phi_{n}^{\epsilon} \to \phi_{\ell}^{\epsilon}$ in $W^{2,2}(\Oe)$. By Sobolev embeddings, we conclude that if $N>4$, then $\phi_{n}^{\epsilon}\to \phi_{\ell}^{\epsilon}$ in $L^{q_{0}}(\Oe)$ with $q_{0}:=\frac{2N}{N-4}>2$. Thus, by a bootstrap argument, we obtain
\[
\phi_{n}^{\epsilon} \to \phi_{\ell}^{\epsilon} \quad \text{in } W^{2,p}(\Oe), \quad \text{for all } p>N.
\]

On the other hand, we also have
$$
    \Lop_2 \phi_n^\epsilon + B_2(x)\phi_n^\epsilon - \lambda_n\phi_n^\epsilon = \epsilon \quad \text{in } \Onl,
$$
or, equivalently,
\begin{align}
    -2\int_{\Onl} J(x-y)\phi_n^\epsilon(y)dy - \int_{\Oe} J(x-y)\phi_n^\epsilon(y)dy + (a(x)-\lambda_{n})\phi_n^\epsilon = \epsilon \quad \text{in } \Onl. \label{2.29}
\end{align}
Due to the compactness of the integral operator, there exists a subsequence (still indexed by $n$) such that
\[
2 \int_{\Onl}J(x-y)\phi_n^{\epsilon}(y)dy \to f \quad \text{in } C(\overline{\Omega}_{n\ell}).
\]
In particular, $f \ge 0$ in $\Onl$. Similarly, using the convergence of $\phi_n^\epsilon|_{\Oe}$, we conclude that
\[
\int_{\Oe}J(x-y)\phi_n^{\epsilon}(y) dy \to \int_{\Oe}J(x-y)\phi_\ell^{\epsilon}(y) dy \quad \text{in } C(\overline{\Omega}_{n\ell}).
\]
Thus, passing to the limit as $n \to +\infty$ in \eqref{2.29}, we obtain
\begin{multline*}
\displaystyle    \phi_n^\epsilon = \frac{\displaystyle \epsilon + 2\int_{\Onl} J(x-y)\phi_n^\epsilon(y)dy + \int_{\Oe} J(x-y)\phi_n^\epsilon(y)dy}{a(x)-\lambda_n} \to \cr\frac{ \displaystyle \epsilon + f+\int_{\Oe}J(x-y)\phi_\ell^{\epsilon}(y) dy}{a(x)-\lambda_p} =: \phi_{n\ell}^\epsilon \quad \text{in } C(\overline{\Omega}_{n\ell}).
\end{multline*}
This directly implies that $\phi_{n\ell}^{\epsilon} \in \operatorname{int}P_{n\ell}$. Thus, denoting $\phi^\epsilon := (\phi_\ell^\epsilon, \phi_{n\ell}^\epsilon) \in X_\mathfrak{B}$ and letting $n \to \infty$ in \eqref{P_ep} yields
\begin{align*}
\begin{cases}
(L+B(x)-\lambda_p)\phi^\epsilon = \epsilon & \text{in } \Omega, \\
\mathfrak{B} \phi^\epsilon = 0 & \text{on } \partial \Oe. 
\end{cases}
\end{align*}
Hence,
\[
\frac{L\phi^\epsilon+B(x)\phi^{\epsilon}}{\phi^\epsilon} = \lambda_p + \frac{\epsilon}{\phi^\epsilon} \quad \text{in } \Omega,
\]
implying that
\[
 \lambda_p = \sup_{\phi\in\mathfrak{P}} \inf_{ \Omega} \frac{(L+B(x))\phi}{\phi}\ge \inf_{\Omega} \frac{L\phi^\epsilon+B(x)\phi^{\epsilon}}{\phi^\epsilon} \ge \lambda_p + \inf_{\Omega} \frac{\epsilon}{\phi^\epsilon} > \lambda_p,
\]
which is a contradiction.

We now analyze \textbf{Case 2}. Let $\epsilon_n \to 0^+$ and define
\[
\phi_n := \phi_{\lambda_{\epsilon_n}}^{\epsilon_n}, \quad \psi_n := \frac{\phi_n}{|\phi_n|_2}, \quad \text{and} \quad c_n := \frac{\epsilon_{n}}{|\phi_n|_2}.
\]
Then $|\psi_n|_2 = 1$ and $c_n \to 0$ as $n \to \infty$. Furthermore, for a subsequence, we have $\lambda_n := \lambda_{\epsilon_n} \to \overline{\lambda} \in [\lambda_p-1, \lambda_p]$.
Dividing equation \eqref{P_ep} by $|\phi_n|_{2}$, we obtain
\begin{align}
\begin{cases}
\Lop_1 \psi_n + B_1(x)\psi_n - \lambda_n \psi_n = c_n & \text{in } \Oe, \\
\mathfrak{B} \psi_n = 0 & \text{on } \partial \Oe, \\
\Lop_2 \psi_n + B_2(x)\psi_n - \lambda_n \psi_n = c_n & \text{in } \Onl. \label{2.31}
\end{cases}
\end{align}
Proceeding as in \textbf{Case 1}, it follows that $\psi_n \to \psi_\ell$ in $W^{2,p}(\Oe)$ for $p>N$.
On the other hand, in $\Onl$ we have
\[
-2\int_{\Onl} J(x-y)\psi_n(y)dy - \int_{\Oe} J(x-y)\psi_n(y)dy + (a(x)-\lambda_n)\psi_n = c_n.
\]
Since $|\psi_n|_{2} = 1$, we can argue as in \textbf{Case 1} to conclude that
\[
2 \int_{\Onl}J(x-y)\psi_n(y)dy \to g \quad \text{in } C(\overline{\Omega}_{n\ell})
\]
and
$$
    \psi_n \to \frac{ g+\int_{\Oe}J(x-y)\psi_\ell(y) dy}{a(x)-\overline{\lambda}} =: \psi_{n\ell} \quad \text{in } C(\overline{\Omega}_{n\ell}).
$$
Thus $\psi_{n\ell} \in \operatorname{int}P_{n\ell}$. Defining $\psi := (\psi_\ell, \psi_{n\ell}) \in X_\mathfrak{B}$, we have $\psi \neq 0$, and taking the limit as $n \to +\infty$ in \eqref{2.31}, we conclude that
\[
    \begin{cases}
(L+B(x)-\overline{\lambda})\psi = 0 & \text{in } \Omega, \\
\mathfrak{B} \psi = 0 & \text{on } \partial \Oe. 
\end{cases}
\]
Therefore, $\psi$ is a positive supersolution associated with $(L+B-\overline{\lambda},\mathfrak{B},\Omega)$. Following the proof of Theorem \ref{SS-SMP}, it follows that $\psi \gg 0$. Finally, if $\overline{\lambda} < \lambda_p$, we would have $\operatorname{ker}[L+B-\overline{\lambda}] = \{0\}$, implying $\psi = 0$, a contradiction. Consequently, $\overline{\lambda} = \lambda_p$, and $\varphi_p = \psi \gg 0$ is the desired eigenfunction.

We now prove item (ii). By Proposition \ref{Prop:lambda-av-principal}, if there exists another principal eigenvalue $\lambda \in \mathbb{R}$ of $(L+B,\mathfrak{B},\Omega)$, then necessarily $\lambda_{p} \leq \lambda < m$. Suppose, by contradiction, that $\lambda > \lambda_{p}$. Thus, there exists an eigenfunction $\phi_{\lambda} \in X_{\mathfrak{B}}$ with $\phi_{\lambda} > 0$ satisfying 
    \begin{align*}
        \begin{cases}
            (L+B(x))\phi_{\lambda} = \lambda\phi_{\lambda} > \lambda_{p} \phi_{\lambda}  &\text{in } \Omega, \\
            \mathfrak{B}\phi_{\lambda} = 0  &\text{on } \partial \Oe.
        \end{cases}
    \end{align*}
Therefore, $\phi_{\lambda}$ is a positive strict supersolution for $(L+B-\lambda_{p}, \mathfrak{B}, \Omega)$. By Corollary \ref{SS-SMP}, $(L+B-\lambda_{p}, \mathfrak{B}, \Omega)$ satisfies the strong maximum principle. In particular, $-\phi_{p} < 0$ satisfies
     \begin{align*}
        \begin{cases}
           (L+B(x)-\lambda_{p})(-\phi_{p}) = 0 &\text{in } \Omega, \\
           \mathfrak{B}(-\phi_{p}) = 0  &\text{on } \partial\Oe,
        \end{cases}
    \end{align*}
meaning that $-\phi_{p}$ is also a supersolution for the operator $(L+B-\lambda_{p}, \mathfrak{B}, \Omega)$, which implies $-\phi_{p} \gg 0$, a contradiction. This completes the proof of the uniqueness of the principal eigenvalue.

To conclude, we will prove item (iii). Since $\lambda_p < m$, it follows that $L+B - \lambda_p$ is a Fredholm operator of index zero. Furthermore, because $L+B-\lambda$ is bijective for all $\lambda < \lambda_p$, we have $\lambda \in \rho(L+B)$.

On the other hand, since Fredholm operators of index zero are also semi-Fredholm operators, it follows from \cite[Theorem 5.31]{pub.1016619032} that there exists an $\epsilon > 0$ such that 
\[
\dim \ker[L+B-\lambda] = \text{constant}, \quad \text{for all } \lambda \in \mathbb{C} \text{ with } 0<|\lambda-\lambda_{p}|<\epsilon.
\]

Consequently, $\dim \ker [L+B-\lambda]=0$ for all $\lambda \in \mathbb{C}$ such that $0<|\lambda-\lambda_{p}|<\epsilon$. Since the set of Fredholm operators of index zero is open in the operator topology, it follows that $L+B-\lambda$ remains a Fredholm operator of index zero for all $|\lambda-\lambda_{p}|<\epsilon$. Thus, $\lambda \in \rho(L+B)$ for all $0<|\lambda-\lambda_{p}|< \epsilon$.

Next, we ensure that  $\lambda_{p}$ is simple. First, we show that $\dim \ker[L+B-\lambda_{p}] = 1$. Let $\varphi \in \ker[L+B-\lambda_{p}]$ with $\varphi \neq 0$. Without loss of generality, we may assume $\varphi < 0$ in somewhere. By (i), $\varphi_{p} \gg 0$ is a positive supersolution for the operator $(L+B-\lambda_{p}, \mathfrak{B}, \Omega)$. Since $\varphi$ is also a supersolution, Theorem \ref{Teorema C} implies the existence of $\tau < 0$ such that $\varphi = \tau \varphi_{p}$. Hence, $\ker[L+B-\lambda_{p}] = \operatorname{span}\{\varphi_p\}$.

Assume now that $u \in X_{\mathfrak{B}}$ is a solution of \eqref{simple-eigenvalue}. Multiplying \eqref{simple-eigenvalue} by $\varphi_p\gg0$, integrating over $\Omega$ and using integration by parts combined with the symmetry of $J$ we obtain
$$
0<\langle\varphi_p,\varphi_p\rangle_2=\langle (L+B(x)-\lambda_p)u,\varphi_p\rangle_2 = \langle u,(L+B(x)-\lambda_p)\varphi_p\rangle_2 =0,
$$
which is a contradiction. This completes the proof.

\end{proof}

\begin{remark} {\rm It is worth noting that Theorem \ref{Th: caso l_p<m} (iii) holds true even without a symmetry assumption on $J$. The proof of this broader case requires the use of the adjoint operator of $L$ and some additional properties (see, e.g.,  \cite[Theorem 2.3]{MR3285831}). For the sake of clarity and simplicity, we restrict our exposition to the symmetric case.}
\end{remark}

The following result is a version of the characterization of the maximum principle from \cite[Theorem 2.5]{LG-JDE96} (see also \cite{bookSMP}) and \cite[Theorem 2.3]{MR3583503}.

\begin{theorem}\label{Caracte PM}
    Assume that $\lambda_{p} < m$. The following statements are equivalent:
    \begin{enumerate}
    \item[(i)] $\lambda_{p} > 0$.
    \item[(ii)] $(L+B, \mathfrak{B}, \Omega)$ admits a positive strict supersolution.
    \item[(iii)] $(L+B, \mathfrak{B}, \Omega)$ satisfies the Strong Maximum Principle.
\end{enumerate}
\end{theorem}

\begin{proof}
    To see that (i) implies (ii), it suffices to observe that, since $\lambda_{p} < m$, an eigenfunction $\varphi_{p} \gg 0$ associated with $\lambda_p > 0$ is a positive strict supersolution for the operator $(L+B, \mathfrak{B}, \Omega)$.

    That (ii) implies (iii) follows directly from Corollary \ref{SS-SMP}. 
     
    Finally, we prove that (iii) implies (i). Suppose, by contradiction, that $\lambda_{p} \le 0$ and let $\varphi_p \gg 0$ be the eigenfunction associated with $\lambda_p$. Thus, we have that $-\varphi_{p} < 0$ satisfies
    \begin{align*}
         \begin{cases}
             (L+B(x))(-\varphi_{p}) = -\lambda_{p}\varphi_{p} \ge 0 \ \ &\text{in } \Omega, \\ 
             \mathfrak{B}(-\varphi_{p} ) \ge 0 \ \ &\text{on } \partial \Oe,
         \end{cases}
    \end{align*}
    which means that $-\varphi_{p}$ is a supersolution for the operator $(L+B, \mathfrak{B}, \Omega)$. It follows from the strong maximum principle that $-\varphi_{p} \gg 0$, which is a contradiction. This completes the proof.
\end{proof}

We conclude this section by providing a sufficient condition for $\lambda_p < m$. From now on, we will denote
$$M:= \displaystyle\max_{\overline{\Omega}_{n\ell}}a(x). $$
Recall that \[
a(x) := 2 \int_{\Onl} J(x-y)dy + \int_{\Oe} J(x-y)dy + B_2(x) .
\]

\begin{theorem}\label{Teo 3.6}
If there exists a measurable set $A \subset \Omega_{n\ell}$ such that
\begin{align}
   \frac{1}{2} < \inf_{y \in A} \int_A \frac{J(x-y)}{a(x)-m} dx \le \infty \label{2.34}
\end{align}
or if
\begin{align}
    M-m< \inf_{x \in \Onl}\int_{\Onl}J(x-y)dy \label{2.35}
\end{align}
then  $\lambda_{p}<m$.
\end{theorem}
\begin{proof}
Recall that, by Remark \ref{Remark lambdap menor m}, we have $\lambda_p \le m$. Suppose, by contradiction, that $\lambda_{p}=m$. Then, for every $\lambda < m$, there exists a strict positive supersolution $\phi \in \mathfrak{P}$ of $(L+B-\lambda,\mathfrak{B},\Omega )$ (see Remark \ref{Remark-l<lp}).
In particular,  
\[
\Lop_2 \phi + B_2(x)\phi - \lambda_n \phi = -2\int_{\Onl} J(x-y)\phi(y)dy - \int_{\Oe} J(x-y)\phi(y)dy + (a(x)-\lambda)\phi > 0
\]
in $\Onl$.
Since $a(x) \geq m > \lambda$ for all $x \in \Onl$, it follows that
\begin{align}
     \phi(x) > \frac{\displaystyle 2\int_{\Onl} J(x-y)\phi(y)dy + \int_{\Oe}J(x-y)\phi(y)dy}{a(x)-\lambda}  \nonumber \\
     > \frac{ \displaystyle2\int_{\Onl} J(x-y)\phi(y)dy}{a(x)-\lambda} \quad \text{in } \Onl. \label{2.36}
\end{align}
Given a measurable subset $D \subset \Onl$, integrating \eqref{2.36} over $D$ and applying Fubini's Theorem, we obtain
\begin{align*}
\int_D \phi(x) dx &> \int_D \frac{2\int_{\Onl} J(x-y)\phi(y)dy}{a(x)-\lambda} dx \\ 
& \ge 2\int_D \int_D \frac{J(x-y)}{a(x)-\lambda} \phi(y) dy dx \\
&\ge 2 \int_D \left[ \int_D \frac{J(x-y)}{a(x)-\lambda} dx \right] \phi(y) dy \\
&\ge 2\left[ \inf_{y \in D} \int_D \frac{J(x-y)}{a(x)-\lambda} dx \right] \int_D \phi(y) dy.
\end{align*}
As $\int_D \phi(y)dy > 0$, we get
\begin{align}
  \frac{1}{2} > \inf_{y \in D} \int_D \frac{J(x-y)}{a(x)-\lambda} dx \quad \forall \lambda<m. \label{2.37}  
\end{align} 
Choosing $D$ to be a ball with diameter less than $\delta>0$ (see hypothesis \ref{J1}), we have $\inf_D J>0$. Thus, we can infer from \eqref{2.37} that 
$$
\frac{1}{2}> \inf_{z\in D} J(z) \int_D\frac{1}{a(x)-\lambda}dx.
$$
Taking a sequence $\lambda_n \uparrow m$ and using Fatou's Lemma, we find that $1/(a(x)-m) \in L^1(D)$. Furthermore, because $\Onl$ is bounded, we deduce that $1/(a(x)-m) \in L^1(\Omega_{n\ell})$.

Setting $D=A$ in \eqref{2.37}, the continuity of $J$ guarantees the existence of ${y_{\lambda}} \in \overline{A}$ such that
\begin{align*}
    \frac{1}{2} > \int_A \frac{J(x-y_\lambda)}{a(x)-\lambda} dx.
\end{align*}
Let $\{\lambda_n\}_{n\geq1}$ be a sequence such that $\lambda_n \uparrow m$. Without loss of generality, we may assume that $y_{\lambda_n} \to \overline y \in \overline{A}$. Since $1/(a-m) \in L^1(\Onl)$, Lebesgue's Dominated Convergence Theorem yields
\[
\frac{1}{2} \ge \int_A \frac{J(x-\overline{y})}{a(x)-m} dx \ge \inf_{y \in A} \int_A \frac{J(x-y)}{a(x)-m} dx,
\]
which contradicts \eqref{2.34}. This completes the proof.

 Now, in \eqref{2.36}, we choose $x_\lambda \in \overline{\Omega}_{n\ell}$ such that $\phi(x_{\lambda})= \min_{\overline{\Omega}_{n\ell}}\phi$. Hence, we have
 $$\phi(x_\lambda) > \frac{\displaystyle\int_{\Onl} J(x_\lambda-y)\phi(y)dy}{a(x_\lambda)-\lambda_n} \ge 
 \frac{\displaystyle \int_{\Onl} J(x_\lambda-y)dy}{a(x_\lambda)-\lambda_n} \phi(x_\lambda),$$
 implying that
 \begin{align}
 a(x_\lambda) - \lambda_n > \int_{\Onl} J(x_\lambda-y)dy \label{2.42} .
 \end{align}
 Passing to a subsequence if necessary, we may assume that $x_\lambda \to \overline{x} \in \overline{\Omega}_{n\ell}$. Thus, letting $\lambda_n \uparrow m$ in \eqref{2.42}, we obtain that:$$M - m \ge \int_{\Onl} J(\overline{x}-y)dy$$contradicting \eqref{2.35}. This concludes the proof.
\end{proof}

\section{Properties of \texorpdfstring{$\lambda_p$}{lp}} \label{sect-proper}

Inspired by \cite{BNV, BCV}, in this section we establish several useful characterizations and monotonicity properties of $\lambda_p$. To this end, we shall adopt the following notation:
$$\lambda_{p}^{\Omega}[L+B]:=\lambda_{p}.$$ We also introduce the set
\begin{equation}
    \mathcal{A} := \Big\{ \lambda \in \mathbb{R} : \exists \, \phi \in X, \, \phi(x) > 0, \, x \in \Omega,
     \text{ supersolution of } (L + B - \lambda, \mathfrak{B}, \Omega) \Big\} \label{Con A}
\end{equation}
and define
\begin{equation*}
    \mu_p := \mu_{p}^{\Omega}[L+B] = \sup \mathcal{A}.
\end{equation*}

The first result of this section establishes the equivalence between $\mu_{p}^{\Omega}[L+B]$ and $\lambda_{p}^{\Omega}[L+B]$.

\begin{proposition}\label{propriedades}
The equality $$\lambda_p = \mu_{p}$$ holds.
\end{proposition}

\begin{proof}
Let $\lambda \in \mathcal{A}$. Then there exists $\phi \in X$, a supersolution of $(L+B-\lambda, \mathfrak{B}, \Omega)$, such that $\phi(x) > 0$ for all $x \in \Omega$. Thus, it follows that
\[
\lambda \le \inf_{x \in \Omega} \frac{(L+B(x))\phi(x)}{\phi(x)} \le \sup_{\mathfrak{P}} \inf_{x \in \Omega} \frac{(L+B(x))\phi(x)}{\phi(x)} = \lambda_p, \quad \forall \lambda \in \mathcal{A},
\]
which implies that $\mu_p \leq \lambda_p$. Conversely, let $\phi \in \mathfrak{P}$ and define
$$
\overline{\mu} := \inf_{x \in \Omega} \frac{(L+B(x))\phi(x)}{\phi(x)}.
$$
If $\overline{\mu} \in \mathbb{R}$, then $(L+B-\overline{\mu})\phi \geq 0$ in $\Omega$; hence, $\overline{\mu} \in \mathcal{A}$. Thus,
$$
\inf_{x \in \Omega} \frac{(L+B(x)) \phi(x)}{\phi(x)} \leq \mu_p \quad \forall \phi \in \mathfrak{P}.
$$
Consequently, we have $\lambda_p \leq \mu_p$, which completes the proof.
\end{proof}

As a consequence, we obtain several monotonicity properties of $\lambda_p$.

\begin{proposition}\label{Prop 4.2}
The following properties hold:
\begin{enumerate}
    \item[(i)] If $c \in \mathbb{R}$, then $$\lambda_{p}^{\Omega}[L+B+c] = \lambda_{p}^{\Omega}[L+B] + c.$$
   
    \item[(ii)] Let $B_1, D_1 \in L^\infty(\Oe)$ and $B_2, D_2 \in C(\overline{\Omega}_{n\ell})$, and define $B := (B_{1}, B_{2})$ and $D := (D_{1}, D_{2})$. If $B \le D$ in $\Omega$, then
    $$\lambda_{p}^{\Omega}[L+B] \le \lambda_{p}^{\Omega}[L+D].$$

    \item[(iii)] Let $B_{1,n}, B_1 \in L^\infty(\Omega_\ell)$ and $B_{2,n}, B_2 \in C(\overline{\Omega}_{n\ell})$, and define $B_n := (B_{1,n}, B_{2,n})$ and $B := (B_1, B_2)$. If $B_{n} \to B$ in $L^\infty(\Omega_\ell) \times C(\overline{\Omega}_{n\ell})$, then $$\lambda_{p}^{\Omega}[L+B_{n}] \to \lambda_{p}^{\Omega}[L+B].$$
\end{enumerate}
\end{proposition}

\begin{proof}
To prove (i), we note that
\begin{multline*}
    \mathcal{A}_c := \{c + \lambda : \lambda \in \mathcal{A}\} \\
    = \{ \mu : \exists \, \phi \in X, \, \phi(x) > 0, \, x \in \Omega, \text{ supersolution of } (L + B + c - \mu, \mathfrak{B}, \Omega) \},
\end{multline*}
where $\mathcal{A}$ is defined by \eqref{Con A}. Thus,
$$ \mu_{p}^{\Omega}[L+B+c] = \sup \mathcal{A}_c = \sup \mathcal{A} + c = \mu_{p}^{\Omega}[L+B] + c. $$
Since $\mu_p = \lambda_p$, the result follows.

To prove (ii), it suffices to define the sets
$$ \mathcal{A}(L+B) := \{ \lambda \in \mathbb{R} : \exists \, \phi \in X, \, \phi(x) > 0, \, x \in \Omega, \text{ supersolution of } (L + B - \lambda, \mathfrak{B}, \Omega) \} $$
and
$$ \mathcal{A}(L+D) := \{ \lambda \in \mathbb{R} : \exists \, \phi \in X, \, \phi(x) > 0, \, x \in \Omega, \text{ supersolution of } (L + D - \lambda, \mathfrak{B}, \Omega) \} $$
and observe that $\mathcal{A}(L+B) \subset \mathcal{A}(L+D)$.

Finally, we prove (iii). Since $B_{n} \to B$ in $L^\infty(\Omega_\ell) \times C(\overline{\Omega}_{n\ell})$, for a given $\epsilon > 0$, there exists $n_{0} \in \mathbb{N}$ such that for all $n \ge n_{0}$, we have
$$ B(x) - \epsilon \le B_{n}(x) \le B(x) + \epsilon \quad \mbox{a. e. in }\Omega. $$
Applying items (i) and (ii), we obtain
$$ \lambda_{p}^{\Omega}[L+B] - \epsilon \le \lambda_{p}^{\Omega}[L+B_{n}] \le \lambda_{p}^{\Omega}[L+B] + \epsilon, $$
whenever $n \geq n_0$. This implies that $\lambda_{p}^{\Omega}[L+B_{n}] \to \lambda_{p}^{\Omega}[L+B]$, concluding the proof.
\end{proof}
The next result provides a monotonicity property of $\lambda_p$ with respect to the domain. For this purpose, we need to introduce some notation. 

\begin{figure}
\centering\input{fig2 }\caption{Admissible configurations of the subdomains $\Omega_{\ell*}$ and $\Omega_{n\ell*}$.}\label{fig: subDominio}\end{figure}
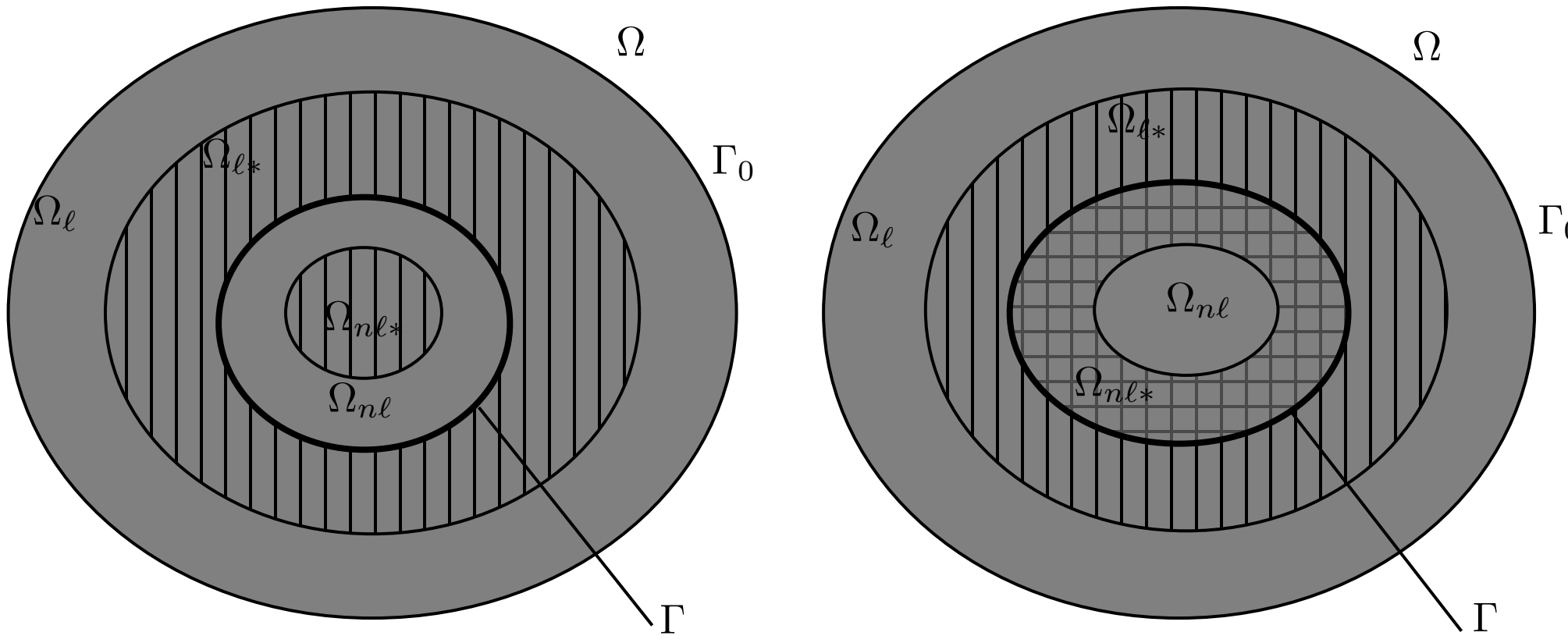

Let $\Omega_{n\ell*} \subset \Onl$ and $\Omega_{\ell*} \subset \Oe$ be regular subdomains with $\operatorname{dist}(\Omega_{\ell*},\Omega_{n\ell*})<\delta$, such that either $\Gamma \subset \partial \Omega_{\ell*}$ or $\Gamma \cap \partial \Omega_{\ell*} = \emptyset$. (see Figure \ref{fig: subDominio}). Setting $\Omega_*= \operatorname{int}(\overline{\Omega_{\ell*} \cup \Omega_{n\ell*}}) \subset \Omega$, we define the sets
\[
    X_*= \Big\{u \in L^2(\Omega_*) : u|_{\Omega_{\ell*}} \in W^{2,p}(\Omega_{\ell*}),\, p > N, \text{ and }  u|_{\overline{\Omega}_{n\ell*}} \in 
    C(\overline{\Omega}_{n\ell*}) \Big\},
\]
$$\mathfrak{P}_*= \Big\{ u \in X_*:  u(x)>0 \mbox{ for all }x \in \Omega_{*} \mbox{ and } \mathfrak{B}_*u \geq 0 \mbox{ on }\partial\Omega_{\ell*}\Big\}$$
and the operators
\[
 L_*u:= \begin{cases} -\Delta u - \displaystyle\int_{\Omega_{n\ell*}}  J(x-y)u(y)dy  + \mathfrak{J}_{n\ell}(x)u  & \text{in } \Omega_{\ell*}, \\ -\displaystyle2\int_{\Omega_{n\ell*}} J(x-y)u(y)dy - \int_{\Omega_{\ell*}} J(x-y)u(y)dy + (2\mathfrak{J}_{n\ell}(x) + \mathfrak{J}_\ell(x))u  & \text{in } \Omega_{n \ell*}, \end{cases}
\]
and
\begin{equation*}
\mathfrak{B}_* u := \begin{cases} 
\displaystyle \frac{\partial u}{\partial \eta} & \text{on } \Gamma,\\
u & \text{on } \partial \Omega_{\ell*} \setminus \Gamma.
\end{cases}
\end{equation*}

\begin{corollary}\label{monotonia-dominio1}
    Under the above conditions, it holds that
    $$
\lambda_p^\Omega[L+B] \leq \lambda_p^{\Omega_*}[L_*+B]. 
    $$
\end{corollary}
\begin{proof}
    Defining 
$$\mathcal{A}_*:= \Big\{\lambda \in \mathbb{R}: \exists  \phi \in X_*, \, \phi(x) > 0,x \in \Omega_*,\ 
\mbox{supersolution of } (L_*+ B-\lambda,\mathfrak{B}_*,\Omega_*) \Big\},$$
given $\lambda \in \mathcal{A}$, there exists a supersolution $\phi \in X$ of $(L+B-\lambda,\mathfrak{B},\Omega)$ such that $\phi(x)>0$ for all $x \in \Omega$. Thus, it is clear that $\phi|_{\Omega_*}$ is a supersolution of $(L_*+B-\lambda,\mathfrak{B}_*,\Omega_*)$. Consequently, $\mathcal{A} \subset \mathcal{A}_*$, which completes the proof.
\end{proof}

We now introduce the following additional definitions:
\begin{align*}
    \overline{\lambda}_{p} := \overline{\lambda}_{p}^\Omega[L+B]=\inf_{\substack{\phi \in X_{\mathfrak{B}}\\ \phi>0}} \sup_{x \in \Omega} \frac{(L+B(x))\phi(x)}{\phi(x)},
\end{align*}
\begin{align*}
   \Tilde{\mathcal{A}} := \Big\{ \lambda : \exists \, \phi \in X_{\mathfrak{B}}, \, \phi > 0, \, \text{subsolution of } (L+ B-\lambda, \mathfrak{B}, \Omega) \Big\}, 
\end{align*}
and
\begin{align*}
\overline{\mu}_{p}:=\overline{\mu}_{p}^{\Omega}[L+B] = \inf \Tilde{\mathcal{A}}.
\end{align*}

Analogously to Proposition \ref{Prop 4.2}, we have:

\begin{proposition}\label{Propo 4.4}
    The equality $\overline{\lambda}_{p} = \overline{\mu}_{p}$ holds.
\end{proposition}

Consequently, it follows that $\overline{\lambda}_p$ satisfies the same monotonicity properties as $\lambda_p$. Specifically:

\begin{proposition}\label{Prop 4.5}
    The following properties hold:
    \begin{enumerate}
        \item [(i)] If $c \in \mathbb{R}$, then $$\overline{\lambda}_{p}^{\Omega}[L+B+c] = \overline{\lambda}_{p}^{\Omega}[L+B] + c.$$
        
        \item [(ii)] Let $B_1, D_1 \in L^\infty(\Omega_\ell)$ and $B_2, D_2 \in C(\overline{\Omega}_{n\ell})$, and define $B := (B_{1}, B_{2})$ and $D := (D_{1}, D_{2})$. If $B \le D$ in $\Omega$, then
        $$\overline{\lambda}_{p}^{\Omega}[L+B] \le \overline{\lambda}_{p}^{\Omega}[L+D].$$

        \item [(iii)] Let $B_{1,n}, B_1 \in L^\infty(\Omega_\ell)$ and $B_{2,n}, B_2 \in C(\overline{\Omega}_{n\ell})$, and define $B_n := (B_{1,n}, B_{2,n})$ and $B := (B_1, B_2)$. If $B_{n} \to B$ in $L^\infty(\Omega_\ell) \times C(\overline{\Omega}_{n\ell})$, then $$\overline{\lambda}_{p}^{\Omega}[L+B_{n}] \to \overline{\lambda}_{p}^{\Omega}[L+B].$$
    \end{enumerate}
\end{proposition}

Let us conclude this section by showing that $\lambda_{p}^{\Omega}[L+B]$ and $\overline{\lambda}_{p}^{\Omega}[L+B]$ coincide.
 
\begin{proposition}\label{Prop 4.6}
The equality $$\lambda_{p}=\overline{\lambda}_{p}$$ holds.  
\end{proposition}

\begin{proof}
First, we will show that $\lambda_{p} \le \overline{\lambda}_{p}$. Suppose, by contradiction, that $\lambda_{p}> \overline{\lambda}_{p}$. Thus, given $\lambda \in (\overline{\lambda}_{p}, \lambda_{p})$, it follows from the definition of $\overline{\lambda}_{p}$ that there exists $\psi \in X_\mathfrak{B}$ such that $\psi>0$ and 
$$ \sup_{x\in\Omega}\frac{(L+B(x))\psi}{\psi}< \lambda.$$
Consequently, 
$$
    \begin{cases}
    (L+B(x)-\lambda)(-\psi)>0  &\mbox{in }  \Omega, \\
    \mathfrak{B} (-\psi) = 0  &\mbox{on }  \partial \Oe.
\end{cases}$$
On the other hand, since $\lambda<\lambda_{p}$, Remark \ref{Remark-l<lp} implies that $(L+B-\lambda, \mathfrak{B}, \Omega)$ satisfies the Strong Maximum Principle. Consequently, $-\psi \gg 0$, which is a contradiction. 

We now prove that $\lambda_{p} \ge \overline{\lambda}_{p}$. If $\lambda_{p}<m$, Theorem \ref{Th: caso l_p<m} ensures the existence of a strictly positive eigenfunction $\varphi_{p}\gg0$ associated with $\lambda_p$. In particular, 
$\varphi_{p}\in X_\mathfrak{B}$ is a subsolution of $(L+B-\lambda_p,\mathfrak{B},\Omega)$, and therefore, $\lambda_{p} \in \Tilde{\mathcal{A}}$. By Proposition \ref{Propo 4.4}, we have $\lambda_{p} \ge \overline{\mu}_{p}=\overline{\lambda}_{p}$. 

Subsequently, we deal with the case $\lambda_{p}=m$. Let $x_0 \in \overline{\Omega}_{nl}$ be such that $a(x_0) =m=\inf_{\Onl}a(x) $. 
 Thus, there exists a sequence $\{x_n\}_{n \ge 1}$ in $\Onl$ such that
 \begin{align*}
 |m - a(x_n)| \le \frac{1}{n}.
 \end{align*}
 Furthermore, since $a \in C(\overline{\Omega}_{n\ell})$, for each $n \in \mathbb{N}$ there exists $\eta_n > 0$ such that
 \begin{align}
 \sup_{B_{\eta_n}(x_n)}|m - a(x)| \le \frac{2}{n}, \quad \forall \ x \in B_{\eta_n}(x_n) \subset \Onl. \label{eqqq 4.7}
 \end{align}
 Now, let us consider a sequence $\{\epsilon_n\}_{n \ge 1}$ such that $\epsilon_n \to 0$ and $\epsilon_n \le \eta_n/2$. We define the following function
 \begin{align*}
 \chi_n(x) := \chi \left( \frac{\normaind{x - x_n}{\infty}}{\epsilon_n} \right),
 \end{align*}
 where $\chi: \mathbb{R}^N \to [0,1]$ is a continuous function satisfying $\chi(x) = 1$ for all $\vert{}x\vert{} \le 1$ and $\chi(x) = 0$ for all $\vert{}x\vert{} \ge 2$. Let us introduce the sequence
 $$a_n(x) := \inf \{a(x), (m - \sup a)\chi_n(x) + \sup a\}.$$
 Notice that
 $$a_n(x) = \begin{cases}  a(x), & \text{if } x \in \overline{\Omega}_{n\ell} \setminus B_{2\epsilon_n}(x_n), \\  m, & \text{if } x \in \overline{\Omega}_{n\ell} \cap B_{\epsilon_n}(x_n). \end{cases}$$
 By the choice of $\epsilon_n$ and condition \eqref{eqqq 4.7}, it follows that$$\normaind{a - a_n}{\infty} \le \sup_{B_{\eta_n}(x_n)}\vert{}m - a(x)\vert{} \to 0.$$
 Moreover, by construction we have $\min_{\overline{\Omega}_{n\ell}} a_n = m$ and $a_n(x) = m$ in $B_{\epsilon_n}(x_n)$. Consequently, $a_n$ satisfies \eqref{2.34}, which implies that $\lambda_p^\Omega[L + B_n] < m$. Thus, from the previous paragraph we obtain
 $$
 \lambda_p^\Omega[L + B_n] = \overline{\lambda}_p^\Omega[L + B_n].$$Finally, by applying Propositions \ref{Prop 4.2} and \ref{Prop 4.5}, we can let $n \to \infty$ in this last equality to conclude that$$\lambda_p^{\Omega}[L + B] = \overline{\lambda}_p^{\Omega}[L + B].$$This completes the proof.
\end{proof}

% 
%-------------------------------------------------------------------------------------------------------
%
%{\color{blue} Sí nos hace falta, abajo una prueba sin la construcción de las $a_n$.}
%
%
%We define sequences $a_n \in C(\overline{\Omega}_{nl})$ and $\epsilon_n>0$ such that 
%$$
%a_n \to a \quad \mbox{in }C(\overline{\Omega}_{nl}), \quad \epsilon_n \to 0^+, \quad \mbox{and}\quad a_n(x)= m \quad \mbox{in } B_{\epsilon_n}(x_0) \cap \overline{\Omega}_{nl}. 
%$$
%Then, defining $B_n=(B_1,B_{2,n})$, where
%$$
%B_{2,n}(x):= a_n(x) - 2\mathfrak{J}_{n\ell}(x)-\mathfrak{J}_{\ell}(x),  
%$$
%we have
%$$
%B_{2,n} \to B_2 \quad \mbox{in }C(\overline{\Omega}_{nl}).
%$$
%Moreover, in $\Onl$,
%$$
%Lu+B_n(x)u= \Lop_2 u + B_{2,n} u= -2\int_{\Onl} J(x-y)u(y)dy - \int_{\Oe} J(x-y)u(y)dy + a_n(x)u.
%$$
%By construction, we have $\min a_n = m$, and $a_n$ satisfies \eqref{2.34}. Consequently, $\lambda_p[L+B_n]<m$. Therefore, by the previous paragraph,
%$$
%\lambda_p[L+B_n] = \overline{\lambda}_p[L+B_n].
%$$
%Now, by Propositions \ref{Prop 4.2} and \ref{Prop 4.5}, we can let $n\to \infty$ in this last equality to obtain
%$$\lambda_{p}^{\Omega}[L+B]=\overline{\lambda}_{p}^{\Omega}[L+B].$$
%This concludes the proof.
%\\
%
%

The previous results also allow us to relate the positivity of $\lambda_p$ to the solvability of the associated linear problem. More precisely:

\begin{corollary}
    Suppose that $\lambda_{p}[L+B]>0$. Then, for every $f \in Y$, there exists a unique $u \in X_\mathfrak{B}$ solving the linear problem
\begin{equation}\label{PL0}
    \begin{cases}
L u (x)  + B(x)u (x) = f(x) & x\in \Omega, \\
\mathfrak{B} u (x) = 0&x\in \partial \Oe.
\end{cases}
\end{equation}
\end{corollary}
\begin{proof}
Since $0<\lambda_{p}[L+B] \le m$, it follows from Remark \ref{Remark-l<lp} that $(L+B, \mathfrak{B}, \Omega)$ admits a strict supersolution. Moreover, we have that$$\ker [L+B]=\{0\}.$$Furthermore, the proof of Theorem \ref{Th: caso l_p<m} guarantees that $L+B: X_{\mathfrak{B}} \to Y$ is a Fredholm operator of index zero. Therefore, we infer that $L+B: X_{\mathfrak{B}} \to Y$ is a bijection. Thus, given $f \in Y$, there exists a unique $u \in X_{\mathfrak{B}}$ satisfying \eqref{PL0}. This concludes the proof.
\end{proof}

It is worth to notice that the existence of a (weak) solution to \eqref{PL0} was proved in \cite[Theorem 2.3]{RossiSIAM} under the assumptions that $B\equiv 0$ and $J$ is symmetric.

\section{Variational characterization of \texorpdfstring{$\lambda_p$}{lp}} \label{sect-variational}

In this section, we work under the assumption that $J$ is symmetric. Our main goal is to show that, under this assumption, a variational characterization for the eigenvalue  $\lambda_{p}^{\Omega}[L+B]$ can be obtained. To this end, we will consider the following space
$$H= \Big\{ u \in L^2(\Omega); u|_{\Oe} \in H^1_{\Gamma_0}(\Oe)\Big\}.$$
Note that $H$ is a Banach space that can be identified with $H_{\Gamma_0}^1(\Oe)\times L^2(\Onl)$. Thus, we can define the following norm
$$
\|u\|_{H}:=\|u|_{\Oe}\|_{H_{\Gamma_0}^1(\Oe)} + \|u|_{\Onl}\|_{L^2(\Onl)}.
$$
\begin{definition}
    We say that $u\in H$ is a weak solution of \eqref{P_l} if $u|_{\Oe}$ and $u|_{\Onl}$ are weak solutions of
 \[\begin{cases} 
\Lop_1 u + B_1(x)u= \lambda u  & \text{in } \Omega_{\ell}, \\ 
\mathfrak{B}u = 0 & \text{on } \partial \Omega_{\ell},
\end{cases} 
\]   
and
$$
\Lop_2 u + B_2(x)u =\lambda u  \quad \text{a.e. in } \Omega_{n \ell},
$$
respectively.
\end{definition}

We also define the bilinear form $b:H \times H \to \mathbb{R}$ given by
$$
b(u,v) = \int_{\Oe} \nabla u \nabla vdx \\
    +\int_{\Onl}\int_{\Omega}J(x-y)(u(y)-u(x))(v(y)-v(x))
    dydx
    +\int_{\Omega}B(x)u vdx.
$$
With this, we introduce the following definition
\begin{align*}
    \lambda_{v}^{\Omega}:=&\inf_{\substack{\phi \in H\\ \phi\ne0}} \frac{\displaystyle \int_{\Oe} |\nabla \phi |^{2}dx+\int_{\Onl}\int_{\Omega}J(x-y)(\phi(y)-\phi(x))^{2}dy dx+\int_{\Omega}B(x)\phi(x)^{2}dx}{\displaystyle\int_{\Omega}\phi^{2}dx}\\
     = &\inf_{\substack{\phi \in H\\ \phi \ne 0}}\frac{\displaystyle b(\phi,\phi)}{\displaystyle \int_{\Omega}\phi^{2}} = \inf_{\substack{\phi \in H\\ |\phi|_2 =1}}b(\phi,\phi)
\end{align*}

We now present an auxiliary result.
\begin{lemma}\label{Lemma 4.7}
    If $J$ is symmetric, then  
    $$
    \left\langle (L+B(x))\phi,\phi \right\rangle_2 = b(\phi,\phi) \quad \forall \phi \in X_\mathfrak{B}.
    $$
\end{lemma}

\begin{proof}
    
   For each $\phi \in X_\mathfrak{B}$, applying integration by parts, we conclude that
\begin{multline}
     \left\langle (L+B(x))\phi,\phi \right\rangle_2= \int_{\Oe} |\nabla \phi |^{2}dx-\int_{\Oe}\int_{\Onl}J(x-y)(\phi(y)-\phi(x))dy\phi(x)dx \\+\int_{\Oe}B_{1}(x)\phi(x)^{2}dx-2 \int_{\Onl}\int_{\Onl}J(x-y)(\phi(y)-\phi(x))dy\phi(x)dx \\ -\int_{\Onl}\int_{\Oe}J(x-y)(\phi(y)-\phi(x))dy\phi(x)dx+\int_{\Onl}B_{2}(x)\phi(x)^{2}dx. \label{4.9}
\end{multline}
   Using Fubini's theorem and the symmetry of J, we obtain
  \begin{multline}
      -\int_{\Oe}\int_{\Onl}J(x-y)(\phi(y)-\phi(x))\phi(x)dydx-\int_{\Onl}\int_{\Oe}J(x-y)(\phi(y)-\phi(x))\phi(x)dydx \\
      = \int_{\Onl}\int_{\Oe}J(x-y)(\phi(y)-\phi(x))\phi(y)dydx-\int_{\Onl}\int_{\Oe}J(x-y)(\phi(y)-\phi(x)) \phi(x)dydx \\ = 
      \int_{\Onl}\int_{\Oe}J(x-y)(\phi(y)-\phi(x))^{2}dydx   \label{4.11}
  \end{multline}
and
\begin{multline}
           -2 \int_{\Onl}\int_{\Onl}J(x-y)(\phi(y)-\phi(x))dy\phi(x)dx \\ = \int_{\Onl}\int_{\Onl}J(x-y)(\phi(y)-\phi(x))\phi(y)dydx-\int_{\Onl}\int_{\Onl}J(x-y)(\phi(y)-\phi(x))\phi(x)dydx
            \\ = \int_{\Onl}\int_{\Onl}J(x-y)(\phi(y)-\phi(x))^{2}dydx. \label{4.12}
       \end{multline}
Thus, substituting \eqref{4.11} and \eqref{4.12} into \eqref{4.9} we obtain 
\begin{align*}
\langle (L+B)\phi, \phi \rangle_2 
&= \int_{\Oe}|\nabla\phi|^{2}+\int_{\Onl}\int_{\Oe}J(x-y)(\phi(y)-\phi(x))^{2}dydx \\
&+ \int_{\Onl}\int_{\Onl}J(x-y)(\phi(y)-\phi(x))^{2}dydx +\int_{\Omega}B(x)\phi^{2}(x)dx\\
&= \int_{\Oe} |\nabla \phi |^{2}dx+\int_{\Onl}\int_{\Omega}J(x-y)(\phi(y)-\phi(x))^{2}dydx+\int_{\Omega}B(x)\phi^{2}(x)dx \\ 
&= b(\phi,\phi),
\end{align*}
 which completes the proof.
\end{proof}

The next result establishes a variational characterization for $\lambda_{p}$ when $J$ is symmetric.
\begin{theorem}\label{Th-CV}
    If $J$ is symmetric, then the equality $$\lambda_{p}=\lambda_{v}$$ holds.
\end{theorem}
\begin{proof}
By Propositions \ref{Propo 4.4} and \ref{Prop 4.6}, it suffices to show that $$\lambda_v= \overline{\mu}_{p}.$$
First, we will show that $\lambda_v \le \overline{\mu}_{p}$. Given $\lambda>\overline{\mu}_{p}$, there exists $\phi\in X_\mathfrak{B}$ such that $\phi>0$ and
\[\begin{cases}
    (L+B(x))\phi\leq \lambda \phi  &\mbox{in }   \Omega, \\
    \mathfrak{B}\phi=0 &\mbox{on }   \partial\Oe.
\end{cases}\]
Multiplying the inequality above by $\phi>0$, integrating over $\Omega$, and using Lemma \ref{Lemma 4.7}, we obtain
$$\lambda_v\leq\frac{b(\phi,\phi)}{|\phi|_2^2}=\frac{\left\langle(L+B(x))\phi,\phi\right\rangle_2}{|\phi|_2^2} \le \lambda \quad \forall \lambda>\overline{\mu}_p.$$
Therefore, $\lambda_v \le \overline{\mu}_{p}$. 

Now suppose, by contradiction, that $\lambda_v<\overline{\mu}_{p}  = \lambda_p$. It follows from the definition of $\lambda_v$ that there exists a sequence $\{\phi_{n}\}_{n \ge 1}$ in $H$ with $\int_{\Omega}\phi_{n}^{2}=1$ such that  
\begin{align}\label{sequn minimizante}
    b(\phi_n,\phi_n) \to \lambda_v. 
\end{align}
Define the symmetric bilinear map $a:H\times H \to \mathbb{R}$ given by
\begin{multline*}
    a(\phi,v):=b(\phi,v) -\lambda_v\langle \phi,v\rangle_2
    =\int_{\Oe} \nabla \phi \nabla vdx \\
    +\int_{\Onl}\int_{\Omega}J(x-y)(\phi(y)-\phi(x))(v(y)-v(x))
    dydx
    +\int_{\Omega}B(x)\phi vdx-\lambda_v\int_{\Omega}\phi vdx. 
\end{multline*}
The definition of $\lambda_v$ also implies that $a$ is positive semidefinite. Thus, by the Cauchy-Schwarz inequality, 
$$    |a(\phi,v)|^{2} \le  a(\phi,\phi)a(v,v) \quad \forall\phi,v \in H.$$
A direct calculation shows that there exists $\gamma>0$ such that
$$a(v,v) \le \gamma (\normaind{v}{H^{1}(\Oe)}^{2}+\|v\|_{L^{2}(\Onl)}^{2}).$$
 Indeed, note that
\begin{multline*}
   a(v,v)\le\int_{\Oe}|\nabla v|^{2}+\int_{\Onl}\int_{\Omega}J(x-y)(v(y)-v(x))^{2}dydx  \\ + \int_{\Omega}B(x)v^2+|\lambda_{v}|\int_{\Omega}v^{2}(x)dx. 
\end{multline*}
Thus, 
\begin{multline*}
    a(v,v) \le \int_{\Oe}|\nabla v|^{2}+\int_{\Onl}\int_{\Omega}J(x-y)v(y)^{2}dydx+\int_{\Onl}\int_{\Omega}J(x-y)v(x)^{2}dydx \\ + \int_{\Oe}B_{1}(x)v^{2}+\int_{\Onl}v(x)^{2}+|\lambda_{v}|\int_{\Oe}v(x)^{2}dx+|\lambda_{v}|\int_{\Onl}v^{2}(x)dx+\int_{\Onl}B_{2}(x)v^{2}(x)   \\ \le
    |\nabla v|_{L^{2}(\Oe)}^{2}+\normaind{J}{\infty}K_{0}|v|_{L^{2}(\Oe)}+\normaind{B_{1}}{\infty}|v|_{L^{2}(\Onl)}^{2}+ |\lambda_{v}|\int_{\Onl}v^{2}dx
    \\  +K_{1}\normaind{J}{\infty}|v|_{L^{2}(\Onl)}^{2} +  K_{2} |v|_{L^{2}(\Onl)}^{2}+|\lambda_{v}||v|_{L^{2}(\Onl)}^{2}+\normaind{B_2}{\infty}|v|_{L^{2}(\Onl)}^{2}.
\end{multline*}
 Hence, setting $\gamma:=\max\{1,K_{0}\normaind{J}{\infty} , K_{1}\normaind{J}{\infty}, K_{2}, \normaind{B_1}{\infty},  \normaind{B_2}{\infty},|\lambda_{v}|\}>0$, we obtain
 $$a(v,v) \le \gamma(\normaind{v}{H^{1}(\Oe)}^{2}+|v|_{L^{2}(\Onl)}^{2}).$$
Consequently,
 $$|a(\phi_{n},v)|^{2} \le \gamma a(\phi_n, \phi_n)(\normaind{v}{H^{1}(\Oe)}^{2}+\|v\|_{L^{2}(\Onl)}^{2}),$$
and therefore,
\begin{equation} \label{weak-convergence}
  \sup_{\normaind{v}{H}\le 1}|a(\phi_n,v)| \le  \gamma^{1/2} a(\phi_n, \phi_n)^{1/2} \to 0.  
\end{equation}
\begin{comment}
    \begin{equation} a(\phi_n,v)=b(\phi_{n}, v ) - \lambda_v \langle \phi_n,v\rangle_2 \to 0 \quad \forall v \in H.
\end{equation}
\end{comment}
Subsequently, since $\{\phi_n\}$ is bounded in $L^2(\Omega)$, it follows from \eqref{sequn minimizante} that
$\{\phi_{n}\}$ is also bounded in $H^{1}_{\Gamma_0}(\Oe)$. Therefore, up to a subsequence still denoted by $n$, we get
\begin{equation}\label{conv-L2Ol}
\phi_{n}\rightharpoonup \phi_\ell\quad \text{in } H^1_{\Gamma_0}(\Oe) \quad \mbox{and} \quad
 \phi_{n} \to \phi_\ell \quad \text{in } L^2(\Oe).
\end{equation}
On the other hand, for $v=(0,v_{n\ell}) \in H_{\Gamma_0}^1(\Oe) \times L^2(\Onl)$, we have
\begin{multline*}
\int_{\Omega_{n\ell}}\int_{\Omega}J(x-y)(\phi_{n}(y)-\phi_{n}(x))(v(y)-v(x)) \, dy \, dx \\
= \int_{\Omega_{n\ell}}\int_{\Omega_{n\ell}}J(x-y)\phi_{n}(y)v_{n\ell}(y) \, dy \, dx - \int_{\Omega_{n\ell}}\int_{\Omega}J(x-y)\phi_{n}(y)v(x) \, dy \, dx \\
- \int_{\Omega_{n\ell}}\int_{\Omega_{n\ell}}J(x-y)\phi_{n}(x)v_{n\ell}(y) \, dy \, dx + \int_{\Omega_{n\ell}}\int_{\Omega}J(x-y)\phi_{n}(x)v(x) \, dy \, dx \\
= \int_{\Omega_{n\ell}}\phi_{n}(y)v_{n\ell}(y)\int_{\Omega_{n\ell}}J(x-y) \, dx \, dy - \int_{\Omega_{n\ell}}v_{n\ell}(x)\int_{\Omega}J(x-y)\phi_{n}(y) \, dy \, dx \\
- \int_{\Omega_{n\ell}}v_{n\ell}(y)\int_{\Omega_{n\ell}}J(x-y)\phi_{n}(x) \, dx \, dy + \int_{\Omega_{n\ell}}\phi_{n}(x)v_{n\ell}(x)(\mathfrak{J}_{n\ell}(x)+\mathfrak{J}_{\ell}(x)) \, dx \\
= \int_{\Omega_{n\ell}}v_{n\ell}(x)\left(L_0\phi_n(x)+(2\mathfrak{J}_{n\ell}(x)+\mathfrak{J}_{\ell}(x))\phi_n(x)\right) \, dx.
\end{multline*}
Consequently,
$$ a(\phi_n,v) = \int_{\Onl} v_{n\ell}(L_0 \phi_n + (a(x) - \lambda_v)\phi_n) dx \quad \forall v=(0,v_{n\ell}) \in H.$$
Defining $w_n:=L_0\phi_{n}+(a(x) - \lambda_v)\phi_n$, it follows from \eqref{weak-convergence} that
\begin{align*}
\|w_n\|_{L^2(\Onl)} &= 
\sup_{\|v_{n\ell}\|_{L^2(\Onl)}\leq 1} |\langle w_n, v_{n\ell} \rangle_{L^2(\Onl)}| \nonumber\\ 
&= \sup_{\|v_{n\ell}\|_{L^2(\Onl)}\leq1} \left|\int_{\Omega_{n\ell}}\left(L_0\phi_{n}(x)+(a(x ) - \lambda_v)\phi_n(x)\right)v_{n\ell}(x)  \, dx \right| \nonumber \\
&=  \sup_{\|v_{n\ell}\|_{L^2(\Onl)}\leq 1} |a(\phi_n,v)| \nonumber\\
&\leq \sup_{\|v\|_{H}\leq 1} |a(\phi_n,v)|\to 0.
\end{align*}
On the other hand, the boundedness of $\{\phi_{n}\}_{n \ge 1}$ combined with \eqref{conv-L2Ol} yields
\begin{multline*}
    -L_0\phi_n=2 \int_{\Onl}J(x-y)\phi_{n}(y)dy+\int_{\Oe}J(x-y)\phi_{n}(y)dy\\ \to f + \int_{\Oe}J(x-y)\phi_{\ell}(y)dy   \quad \mbox{in }C(\overline{\Omega}_{n\ell}).
\end{multline*}
Furthermore, since $\lambda_v <\lambda_p \leq m$, we have
\begin{equation}\label{conv-L2Onl}
    \phi_n = \frac{w_n -L_0 \phi_n}{a(x)-\lambda_v} \to \frac{f+ \int_{\Oe}J(x-y) \phi_\ell(y)dy}{a(x)-\lambda_v} =: \phi_{n\ell} \quad \mbox{in }L^2(\Omega_{n\ell}).
\end{equation}
Thus, defining $\phi = (\phi_\ell, \phi_{n\ell}) \in H$, we see that
$$
0=\lim\|w_n\|_{L^2(\Onl)} = \|L_0 \phi + (a(x)  -\lambda_v)\phi_{n\ell}\|_{L^2(\Onl)}.
$$
In particular, 
$$
\phi_{n\ell} = \frac{L_0\phi}{a(x)-\lambda_v} \in C(\overline{\Omega}_{n\ell})
$$
and
$$
\Lop_2 \phi + B_2(x) \phi_{n\ell} - \lambda_v \phi_{n\ell}=
L_0 \phi + (a(x)  -\lambda_v)\phi_{n\ell}=0 \quad \mbox{in } \Onl.
$$

Similarly, taking $v=(v_\ell,0) \in H_{\Gamma_0}^1(\Oe)\times L^2(\Onl)$ in \eqref{weak-convergence} and using \eqref{conv-L2Ol} and \eqref{conv-L2Onl}, we find that
\begin{multline*}
   0 =\lim  a(\phi_n,v)\\ = 
    \int_{\Oe}\nabla\phi_\ell \nabla v_\ell- \int_{\Oe}v_\ell(x)  \int_{\Omega_{n\ell}} J(x-y)(\phi_{n\ell}(y)-\phi_\ell(x))dy dx + \int_{\Oe} (B_1(x)-\lambda_v)v_\ell\phi_\ell,
\end{multline*}
which means that $\phi_\ell \in H_{\Gamma_0}^1(\Oe)$ is a weak solution of 
\[\begin{cases}
        -\Delta u-\displaystyle \int_{\Omega_{n\ell}} J(x-y)(\phi_{n\ell}(y)-u(x))dy+B_1(x) u =\lambda_{v} u  &\mbox{in }   \Omega_{\ell}, \\
        \mathfrak{B}u=0 &\mbox{on }   \partial\Oe.
    \end{cases}\]
By elliptic regularity, we have $\phi_\ell \in W^{2,p}(\Oe)$ for $p>N$. Consequently, $\phi=(\phi_\ell,\phi_{n\ell}) \in X_\mathfrak{B}$ is a solution of \eqref{P_l}. Thus,
$$
\phi \in \operatorname{ker}[L+B-\lambda_v].
$$
Since $\lambda_v<\lambda_p$, Remark \ref{Remark-l<lp} implies that $\phi =0$. However, from \eqref{conv-L2Ol} and \eqref{conv-L2Onl} we obtain
$$
1=\lim|\phi_n|_2 = |\phi|_2,
$$
which is a contradiction. This completes the proof.
\end{proof}

\subsection*{Funding}

\noindent W. C. has been  supported by CNPq Proc. 310664/2023-3.  J. P. was financed in part by CAPES Proc. 88881.125193/2025-01 and 88887.328317/2026-00. A. S. has been supported under Grants PID2023-149509NB-I00 and  CEX2024-001517-M - Apoyo a Unidades de Excelencia María de Maeztu  funded by MICIU/AEI. A. S. also has been supported by University of Sevilla under grants SOL2024-31708 and SOL2024-31596. 

%\bibliographystyle{acm} % Estilo de referência simples 
%\bibliography{myref} % Use o exemplo de arquivo bibliográfico sample.bib

%---------------------------------
% BIBLIOGRAFIA
%---------------------------------

\end{document}

%% file: fig1.tex
\tikzset{every picture/.style={line width=0.75pt}} %set default line width to 0.75pt        

\begin{tikzpicture}[x=0.75pt,y=0.75pt,yscale=-1,xscale=1]
%uncomment if require: \path (0,340); %set diagram left start at 0, and has height of 340

%Shape: Ellipse [id:dp8998697591784376] 
\draw  [fill={rgb, 255:red, 128; green, 128; blue, 128 }  ,fill opacity=1 ] (43.67,99) .. controls (43.67,58.68) and (90.23,26) .. (147.67,26) .. controls (205.1,26) and (251.67,58.68) .. (251.67,99) .. controls (251.67,139.32) and (205.1,172) .. (147.67,172) .. controls (90.23,172) and (43.67,139.32) .. (43.67,99) -- cycle ;
%Shape: Ellipse [id:dp5425504589120801] 
\draw  [fill={rgb, 255:red, 155; green, 155; blue, 155 }  ,fill opacity=1 ] (102.67,97.5) .. controls (102.67,78.45) and (126.17,63) .. (155.17,63) .. controls (184.16,63) and (207.67,78.45) .. (207.67,97.5) .. controls (207.67,116.55) and (184.16,132) .. (155.17,132) .. controls (126.17,132) and (102.67,116.55) .. (102.67,97.5) -- cycle ;
%Shape: Ellipse [id:dp9886975774395605] 
\draw  [fill={rgb, 255:red, 128; green, 128; blue, 128 }  ,fill opacity=1 ] (291.67,100) .. controls (291.67,59.68) and (338.23,27) .. (395.67,27) .. controls (453.1,27) and (499.67,59.68) .. (499.67,100) .. controls (499.67,140.32) and (453.1,173) .. (395.67,173) .. controls (338.23,173) and (291.67,140.32) .. (291.67,100) -- cycle ;
%Shape: Ellipse [id:dp3010007038278988] 
\draw  [fill={rgb, 255:red, 155; green, 155; blue, 155 }  ,fill opacity=1 ] (350.67,98.5) .. controls (350.67,79.45) and (374.17,64) .. (403.17,64) .. controls (432.16,64) and (455.67,79.45) .. (455.67,98.5) .. controls (455.67,117.55) and (432.16,133) .. (403.17,133) .. controls (374.17,133) and (350.67,117.55) .. (350.67,98.5) -- cycle ;

%Shape: Polygon Curved [id:ds5606025013776861] 
\draw  [fill={rgb, 255:red, 155; green, 155; blue, 155 }  ,fill opacity=1 ] (244.42,228.02) .. controls (244.91,180.04) and (85.12,204.62) .. (114.09,228.02) .. controls (143.05,251.42) and (154.15,254.93) .. (114.09,298.22) .. controls (74.02,341.51) and (246,344.08) .. (244.42,298.22) .. controls (242.85,252.35) and (243.94,275.99) .. (244.42,228.02) -- cycle ;
%Straight Lines [id:da9521566336063034] 
\draw    (244,259.25) -- (284,259.25) ;
%Shape: Regular Polygon [id:dp7851792042114863] 
\draw  [fill={rgb, 255:red, 155; green, 155; blue, 155 }  ,fill opacity=1 ] (283.04,229.02) .. controls (282.56,181.04) and (442.34,205.62) .. (413.38,229.02) .. controls (384.42,252.42) and (373.31,255.93) .. (413.38,299.22) .. controls (453.45,342.51) and (281.47,345.08) .. (283.04,299.22) .. controls (284.62,253.35) and (283.53,276.99) .. (283.04,229.02) -- cycle ;

% Text Node
\draw (208,21.4) node [anchor=north west][inner sep=0.75pt]    {$\Omega $};
% Text Node
\draw (54,82.4) node [anchor=north west][inner sep=0.75pt]    {$\Omega _{\ell }$};
% Text Node
\draw (149.67,102.4) node [anchor=north west][inner sep=0.75pt]    {$\Omega _{n\ell }$};
% Text Node
\draw (456,22.4) node [anchor=north west][inner sep=0.75pt]    {$\Omega $};
% Text Node
\draw (405.17,101.9) node [anchor=north west][inner sep=0.75pt]    {$\Omega _{\ell }$};
% Text Node
\draw (305.67,87.4) node [anchor=north west][inner sep=0.75pt]    {$\Omega _{n\ell }$};
% Text Node
\draw (246.67,241.57) node [anchor=north west][inner sep=0.75pt]  [font=\fontsize{0.88em}{1.06em}\selectfont]  {$\operatorname{d< \delta }$};
% Text Node
\draw (195,221.4) node [anchor=north west][inner sep=0.75pt]    {$\Omega _{\ell }$};
% Text Node
\draw (292.67,218.4) node [anchor=north west][inner sep=0.75pt]    {$\Omega _{n\ell }$};

\end{tikzpicture}

%% file: fig2.tex
% Pattern Info
 
\tikzset{
pattern size/.store in=\mcSize, 
pattern size = 5pt,
pattern thickness/.store in=\mcThickness, 
pattern thickness = 0.3pt,
pattern radius/.store in=\mcRadius, 
pattern radius = 1pt}
\makeatletter
\pgfutil@ifundefined{pgf@pattern@name@_pscg8tdrm}{
\pgfdeclarepatternformonly[\mcThickness,\mcSize]{_pscg8tdrm}
{\pgfqpoint{-\mcThickness}{-\mcThickness}}
{\pgfpoint{\mcSize}{\mcSize}}
{\pgfpoint{\mcSize}{\mcSize}}
{
\pgfsetcolor{\tikz@pattern@color}
\pgfsetlinewidth{\mcThickness}
\pgfpathmoveto{\pgfpointorigin}
\pgfpathlineto{\pgfpoint{0}{\mcSize}}
\pgfusepath{stroke}
}}
\makeatother

% Pattern Info
 
\tikzset{
pattern size/.store in=\mcSize, 
pattern size = 5pt,
pattern thickness/.store in=\mcThickness, 
pattern thickness = 0.3pt,
pattern radius/.store in=\mcRadius, 
pattern radius = 1pt}
\makeatletter
\pgfutil@ifundefined{pgf@pattern@name@_m353sydl6}{
\pgfdeclarepatternformonly[\mcThickness,\mcSize]{_m353sydl6}
{\pgfqpoint{-\mcThickness}{-\mcThickness}}
{\pgfpoint{\mcSize}{\mcSize}}
{\pgfpoint{\mcSize}{\mcSize}}
{
\pgfsetcolor{\tikz@pattern@color}
\pgfsetlinewidth{\mcThickness}
\pgfpathmoveto{\pgfpointorigin}
\pgfpathlineto{\pgfpoint{0}{\mcSize}}
\pgfusepath{stroke}
}}
\makeatother

% Pattern Info
 
\tikzset{
pattern size/.store in=\mcSize, 
pattern size = 5pt,
pattern thickness/.store in=\mcThickness, 
pattern thickness = 0.3pt,
pattern radius/.store in=\mcRadius, 
pattern radius = 1pt}
\makeatletter
\pgfutil@ifundefined{pgf@pattern@name@_43qek5o6v}{
\pgfdeclarepatternformonly[\mcThickness,\mcSize]{_43qek5o6v}
{\pgfqpoint{-\mcThickness}{-\mcThickness}}
{\pgfpoint{\mcSize}{\mcSize}}
{\pgfpoint{\mcSize}{\mcSize}}
{
\pgfsetcolor{\tikz@pattern@color}
\pgfsetlinewidth{\mcThickness}
\pgfpathmoveto{\pgfpointorigin}
\pgfpathlineto{\pgfpoint{0}{\mcSize}}
\pgfusepath{stroke}
}}
\makeatother

% Pattern Info
 
\tikzset{
pattern size/.store in=\mcSize, 
pattern size = 5pt,
pattern thickness/.store in=\mcThickness, 
pattern thickness = 0.3pt,
pattern radius/.store in=\mcRadius, 
pattern radius = 1pt}\makeatletter
\pgfutil@ifundefined{pgf@pattern@name@_42k3d7w3l}{
\pgfdeclarepatternformonly[\mcThickness,\mcSize]{_42k3d7w3l}
{\pgfqpoint{-\mcThickness}{-\mcThickness}}
{\pgfpoint{\mcSize}{\mcSize}}
{\pgfpoint{\mcSize}{\mcSize}}
{\pgfsetcolor{\tikz@pattern@color}
\pgfsetlinewidth{\mcThickness}
\pgfpathmoveto{\pgfpointorigin}
\pgfpathlineto{\pgfpoint{\mcSize}{0}}
\pgfpathmoveto{\pgfpointorigin}
\pgfpathlineto{\pgfpoint{0}{\mcSize}}
\pgfusepath{stroke}}}
\makeatother
\tikzset{every picture/.style={line width=0.75pt}} %set default line width to 0.75pt        

\begin{tikzpicture}[x=0.75pt,y=0.75pt,yscale=-1,xscale=1]
%uncomment if require: \path (0,274); %set diagram left start at 0, and has height of 274

%Shape: Ellipse [id:dp20929748413039162] 
\draw  [fill={rgb, 255:red, 128; green, 128; blue, 128 }  ,fill opacity=1 ] (38.7,147.62) .. controls (38.7,93.42) and (91.11,49.48) .. (155.76,49.48) .. controls (220.42,49.48) and (272.83,93.42) .. (272.83,147.62) .. controls (272.83,201.82) and (220.42,245.75) .. (155.76,245.75) .. controls (91.11,245.75) and (38.7,201.82) .. (38.7,147.62) -- cycle ;
%Shape: Ellipse [id:dp2005719803560988] 
\draw  [color={rgb, 255:red, 0; green, 0; blue, 0 }  ,draw opacity=1 ][pattern=_pscg8tdrm,pattern size=6pt,pattern thickness=0.75pt,pattern radius=0pt, pattern color={rgb, 255:red, 0; green, 0; blue, 0}][line width=0.75]  (69.91,147.62) .. controls (69.91,108.37) and (108.35,76.55) .. (155.76,76.55) .. controls (203.18,76.55) and (241.61,108.37) .. (241.61,147.62) .. controls (241.61,186.86) and (203.18,218.68) .. (155.76,218.68) .. controls (108.35,218.68) and (69.91,186.86) .. (69.91,147.62) -- cycle ;
%Shape: Ellipse [id:dp40711332971982517] 
\draw  [fill={rgb, 255:red, 128; green, 128; blue, 128 }  ,fill opacity=1 ][line width=1.5]  (106.34,151) .. controls (106.34,128.57) and (127.3,110.39) .. (153.16,110.39) .. controls (179.02,110.39) and (199.99,128.57) .. (199.99,151) .. controls (199.99,173.43) and (179.02,191.61) .. (153.16,191.61) .. controls (127.3,191.61) and (106.34,173.43) .. (106.34,151) -- cycle ;
%Shape: Ellipse [id:dp4255246280961502] 
\draw  [color={rgb, 255:red, 0; green, 0; blue, 0 }  ,draw opacity=1 ][pattern=_m353sydl6,pattern size=6pt,pattern thickness=0.75pt,pattern radius=0pt, pattern color={rgb, 255:red, 0; green, 0; blue, 0}][line width=0.75]  (127.89,147.62) .. controls (127.89,136) and (139.12,126.59) .. (152.98,126.59) .. controls (166.83,126.59) and (178.06,136) .. (178.06,147.62) .. controls (178.06,159.23) and (166.83,168.64) .. (152.98,168.64) .. controls (139.12,168.64) and (127.89,159.23) .. (127.89,147.62) -- cycle ;
%Straight Lines [id:da6494367550072884] 
\draw    (190,178) -- (245.75,248.1) ;
%Shape: Ellipse [id:dp4537683337438344] 
\draw  [fill={rgb, 255:red, 128; green, 128; blue, 128 }  ,fill opacity=1 ] (300.81,147.62) .. controls (300.81,93.42) and (351.99,49.48) .. (415.11,49.48) .. controls (478.24,49.48) and (529.41,93.42) .. (529.41,147.62) .. controls (529.41,201.82) and (478.24,245.75) .. (415.11,245.75) .. controls (351.99,245.75) and (300.81,201.82) .. (300.81,147.62) -- cycle ;
%Shape: Ellipse [id:dp6008739167335049] 
\draw  [color={rgb, 255:red, 0; green, 0; blue, 0 }  ,draw opacity=1 ][pattern=_43qek5o6v,pattern size=6pt,pattern thickness=0.75pt,pattern radius=0pt, pattern color={rgb, 255:red, 0; green, 0; blue, 0}][line width=0.75]  (333.61,146.65) .. controls (333.61,107.4) and (371.13,75.58) .. (417.42,75.58) .. controls (463.72,75.58) and (501.24,107.4) .. (501.24,146.65) .. controls (501.24,185.9) and (463.72,217.71) .. (417.42,217.71) .. controls (371.13,217.71) and (333.61,185.9) .. (333.61,146.65) -- cycle ;
%Straight Lines [id:da912837547574863] 
\draw    (451.48,179.83) -- (477.24,213.01) -- (505.91,249.93) ;
%Shape: Ellipse [id:dp3177948326609692] 
\draw  [pattern=_42k3d7w3l,pattern size=6pt,pattern thickness=0.75pt,pattern radius=0pt, pattern color={rgb, 255:red, 74; green, 74; blue, 74}][line width=1.5]  (360.68,147.62) .. controls (360.68,124.39) and (385.05,105.56) .. (415.11,105.56) .. controls (445.17,105.56) and (469.54,124.39) .. (469.54,147.62) .. controls (469.54,170.84) and (445.17,189.67) .. (415.11,189.67) .. controls (385.05,189.67) and (360.68,170.84) .. (360.68,147.62) -- cycle ;
%Shape: Ellipse [id:dp46035651415837353] 
\draw  [color={rgb, 255:red, 0; green, 0; blue, 0 }  ,draw opacity=1 ][fill={rgb, 255:red, 128; green, 128; blue, 128 }  ,fill opacity=1 ][line width=0.75]  (387.9,146.65) .. controls (387.9,135.03) and (401.12,125.62) .. (417.42,125.62) .. controls (433.73,125.62) and (446.95,135.03) .. (446.95,146.65) .. controls (446.95,158.26) and (433.73,167.68) .. (417.42,167.68) .. controls (401.12,167.68) and (387.9,158.26) .. (387.9,146.65) -- cycle ;

% Text Node
\draw (488.52,55.75) node [anchor=north west][inner sep=0.75pt]    {$\Omega $};
% Text Node
\draw (379.8,163.66) node [anchor=north west][inner sep=0.75pt]  [font=\normalsize,color={rgb, 255:red, 0; green, 0; blue, 0 }  ,opacity=1 ]  {$\Omega _{n\ell *}$};
% Text Node
\draw (409.27,136.57) node [anchor=north west][inner sep=0.75pt]    {$\Omega _{n\ell }$};
% Text Node
\draw (390.3,78.95) node [anchor=north west][inner sep=0.75pt]  [color={rgb, 255:red, 0; green, 0; blue, 0 }  ,opacity=1 ]  {$\Omega _{\ell *}$};
% Text Node
\draw (307.92,113.88) node [anchor=north west][inner sep=0.75pt]    {$\Omega {_{\ell }}$};
% Text Node
\draw (508,239.4) node [anchor=north west][inner sep=0.75pt]    {$\Gamma $};
% Text Node
\draw (529,111.53) node [anchor=north west][inner sep=0.75pt]    {$\Gamma _{0}$};
% Text Node
\draw (232.59,54.28) node [anchor=north west][inner sep=0.75pt]    {$\Omega $};
% Text Node
\draw (138.53,142.56) node [anchor=north west][inner sep=0.75pt]    {$\Omega _{n\ell *}$};
% Text Node
\draw (139.82,168.74) node [anchor=north west][inner sep=0.75pt]    {$\Omega _{n\ell }$};
% Text Node
\draw (99.39,90.31) node [anchor=north west][inner sep=0.75pt]  [color={rgb, 255:red, 0; green, 0; blue, 0 }  ,opacity=1 ]  {$\Omega {_{\ell *}}$};
% Text Node
\draw (44.94,108.59) node [anchor=north west][inner sep=0.75pt]    {$\Omega {_{\ell }}$};
% Text Node
\draw (263.4,92.22) node [anchor=north west][inner sep=0.75pt]    {$\Gamma _{0}$};
% Text Node
\draw (246.68,240.04) node [anchor=north west][inner sep=0.75pt]    {$\Gamma $};

\end{tikzpicture}